\documentstyle[amssymb,12pt]{amsart}

\ifx\mathrm\undefined\let\mathrm\bf\fi
\ifx\mathbf\undefined\let\mathbf\bf\fi

\let\leq\leqslant \let\geq\geqslant
\let\tsize\textstyle \def
\Sum{\sum\limits}

 \let\Dl\Delta
\let\epe\epsilon \let\eps\varepsilon \let\epsilon\eps

\let\la\lambda 
\let\om\omega 
 \let\phi\varphi

\let\si\sigma

\let\zt\zeta

\def\Ti{\widetilde T}
\def\Li{\widetilde L}

\def\Ugh{U_q\widehat{gl(2)}}

\newcommand{\tell}{t_1,\dots,t_l}
\newcommand{\zn}{z_1,\dots,z_n}
\newcommand{\lan}{\la_1,\dots,\la_n}
\newcommand{\Zb}{{{\mathcal Z}}}
\newcommand{\F}{{{\frak F}\,}}

\newcommand{\half}{\frac12}
\newcommand{\Z}{{\Bbb Z}}

\newcommand{\C}{{\Bbb C}}
\newcommand{\Q}{{\Bbb Q}}   
\newcommand{\Ref}[1]{{$($\ref{#1}$)$}}
\newcommand{\bean}{\begin{eqnarray}}
\newcommand{\eean}{\end{eqnarray}}
\newcommand{\be}{\begin{displaymath}}
\newcommand{\ee}{\end{displaymath}}
\newcommand{\bea}{\begin{eqnarray*}}   
\newcommand{\eea}{\end{eqnarray*}}

\newcommand{\Id}{{\operatorname{Id}}}

\newcommand{\res}{{\operatorname{res}}}
\newcommand{\Imag}{{\operatorname{Im}}\,}
\newcommand{\Real}{{\operatorname{Re}}\,}
\newcommand{\Ker}{{\operatorname{Ker}}}

\newcommand{\T}{\otimes}

\newcommand{\vs}{\vspace{1.5\baselineskip}}

\newenvironment{proof}{\noindent{\it Proof\/}:\rm}{$\;\Box$
\par\vs}

\newtheorem
{thm}{Theorem}
\newtheorem{example}[thm]{Example}
\newenvironment{remark}{\noindent{\bf Remark.\/}}{$\;\Box$\par\vs}

\newtheorem
{lemma}[thm]{Lemma}

\newtheorem
{corollary}[thm]{Corollary}

\newcommand{\End}{{\operatorname{End}}}

\begin{document}

\title{Functorial properties of the hypergeometric map} 
\author{E. Mukhin and A. Varchenko}
\maketitle
\bigskip
\centerline{September 1998}

\begin{abstract}
The quantized Knizhnik-Zamolodchikov equation is a difference equation defined in terms of rational $R$ matrices. We describe all singularities of hypergeometric solutions to the qKZ equations.
\end{abstract}

\section{Introduction}
The quantized Knizhnik-Zamolodchikov equation (qKZ) is a system of difference equations. The qKZ equation was introduced in \cite{FR} as an equation for matrix elements of vertex operators of a quantum affine algebra.
A special case of the qKZ equation had been introduced earlier in \cite{S} as equations for form factors in integrable quantum field theory.
Later, the qKZ equation was derived as an equation for correlation functions
in lattice integrable models, cf. \cite{JM} and references therein.

In this paper we consider the rational qKZ equation associated with $sl(2)$.
The qKZ equation with values in a tensor product of $sl(2)$ Verma modules $V(\la_1)\T\dots\T V(\la_n)$ was solved in \cite{TV}, cf. \cite{M}, \cite{R}. The solutions $\Psi(z,\la)$, $z=(z_1,\dots,z_n)$, $\la=(\la_1,\dots,\la_n)$ are meromorphic functions written in terms of hypergeometric integrals, see \Ref{int specific}. Here $V(\la_i)$ is the Verma module with highest weight $\la_i\in\C$. We endow the module $V(\la_i)$ with an evaluation Yangian module structure with complex evaluation parameter $z_i$ and denote $V(z_i,\la_i)$. Set $V(z,\la)=V(z_1,\la_1)\T\dots\T V(z_n,\la_n)$. 
The space of the hypergeometric solutions can be naturally identified with the corresponding tensor product of Verma modules $V_q(\la_1)\T\dots\T V_q(\la_n)$ over the quantum group $U_qsl(2)$, where $q=\exp(\pi i/p)$ and $p$ is the step of the qKZ equation. We endow the module $V_q(\la_i)$ with an evaluation structure of affine quantum group $\Ugh$ module with complex evaluation parameter $z_i$ and denote $V_q(z_i,\la_i)$. Set $V_q(z,\la)=V_q(z_1,\la_1)\T\dots\T V_q(z_n,\la_n)$.

Thus, the hypergeometric solutions define a (hypergeometric) map
\be
{\rm qKZ}(z,\la;p):\;  V_q(z,\la)\to V(z,\la),
\ee
see \Ref{Psi}. 

We call parameters $p,z,\la$ generic if the $\Ugh$ module $V_q(z,\la)$ is irreducible. 
For generic values of parameters, the hypergeometric map is an isomorphism of vector spaces. 


The Yangian module $V(z,\la)$ is reducible iff $z_a-z_b+\la_a+\la_b\in\Z_{\geq 0}$ for some $a,b\in\{1,\dots,n\}$. We describe the properties of the hypergeometric map for these values of parameters.

Case 1.
Let $2\la_1=k\in\Z_{\geq 0}$. The module $V(z,\la)$ has a submodule isomorphic to $V(z,\la^\prime)$, where $\la^\prime=(-\la_1-1,\la_2,\dots,\la_n)$. The module $V_q(z,\la)$ has a submodule isomorphic to $V_q(z,\la^\prime)$. 
Then the qKZ map is still a well defined isomorphism of vector spaces. Moreover, it maps the $\Ugh$ submodule $V_q(z,\la^\prime)$ onto the Yangian submodule $V(z,\la^\prime)$ and the restriction of the map qKZ$(z,\la)$ to $V_q(z,\la^\prime)$ coincides up to a scalar factor with the hypergeometric map qKZ$(z,\la^\prime)$, see Theorem \ref{la1=k}.

Case 2. Let $z_2-z_1+\la_1+\la_2=k\in\Z_{\geq 0}$. The module $V_q(z,\la)$ has a submodule isomorphic to $V_q(z^\prime,\la^\prime)$, where $z^\prime,\la^\prime\in\C^n$ are some new values of parameters. The module $V(z,\la)$ has a nontrivial submodule such that the factor module is isomorphic to $V(z^\prime,\la^\prime)$. Then, the qKZ map is a well defined linear map. The kernel of the qKZ map coincides with the submodule $V_q(z^\prime,\la^\prime)$ and the image of the qKZ map coincides with the proper submodule in $V(z,\la)$, see
Corollary \ref{qKZ 1 to 2}.

Case 3.
Let $z_1-z_2+\la_1+\la_2=k\in\Z_{\geq 0}$. The module $V(z,\la)$ has a submodule isomorphic to $V(z^\prime,\la^\prime)$, where $z^\prime,\la^\prime\in\C^n$ are some new values of parameters. The module $V_q(z,\la)$ has a nontrivial submodule such that the factor module is isomorphic to $V_q(z^\prime,\la^\prime)$. Then the qKZ map has a simple pole at the hyperplane $z_1-z_2+\la_1+\la_2=k$. Let Res be the corresponding residue of the qKZ map. The kernel of Res coincides with the proper submodule of $V_q(z,\la)$ and the image of Res coincides with the proper submodule of $V(z,\la)$. Thus, Res is a linear isomorphism of the factor module
$V_q(z^\prime,\la^\prime)$ to the submodule $V(z^\prime,\la^\prime)$. We prove that the map Res coincides with the hypergeometric map qKZ$(z^\prime,\la^\prime)$ up to a scalar factor, see Example \ref{z1-z2=k}.

It is well known that the intertwinings of Yangian modules and of modules over $\Ugh$ possess similar properties, see Lemmas \ref{same kind}-\ref{2 to 1} and Lemmas \ref{q same kind}-\ref{q 2 to 1}.

Let $V(\tilde{z},\tilde{\la})$ be a Yangian module obtained from the tensor product $V(z,\la)$ by a permutation of factors. It is well known that $V(\tilde{z},\tilde{\la})\simeq V(z,\la)$,  $V_q(\tilde{z},\tilde{\la})\simeq V_q(z,\la)$ and the qKZ map intertwines these isomorphisms
.

Let $z,\la\in\C^n$. The values of parameters $z^\prime,\la^\prime\in\C^n$ such that $V(z,\la)\simeq  V(z^\prime,\la^\prime)$ are parametrized by pairs of permutations $\si,\si^\prime\in{\Bbb S}^n$. We show that the qKZ map intertwines all these isomorphisms as well, see the diagram in Theorem \ref{extended monodromy}. This observation allows us to reduce cases 2 and 3 to case 1.

We apply the above observation to include the qKZ equation in a bigger compatible system of difference equations which we call the extended qKZ equation. The shifts of the extended qKZ equation generate a group acting in the space $\C^{2n}$ of parameters $(z,\la)$ isomorphic to $\Z^{2n-1}$.

We use the above study to describe all singularities of the hypergeometric solutions, see Theorem \ref{singularities}, and the Remark after it.

In this paper we treat the case $|\kappa|>1$, where $\kappa$ is a parameter of the qKZ equation \Ref{qKZ operators}.
Our results can be carried in the same fashion in the case $\kappa=1$. We also obtain similar results in the case of tensor products of finite dimensional representations, see Remarks after Theorem \ref{singularities}.

The paper is organized as follows. 
We study the functional model of the Yangian module $V(z,\la)$ in Section \ref{rational section}. In Sections \ref{sl(2) section}-\ref{action section} we fix our notations and recall some technical facts from \cite{TV}. In Sections \ref{bases section}-\ref{factor section} we describe the maps of Yangian modules in terms of spaces of functions. We study the functional model of the $\Ugh$ module $V_q(z,\la)$ in Section \ref{trigonometric section}. Section \ref{trigonometric section} is constructed similarly to Section \ref{rational section}. In Sections \ref{qKZ section}-\ref{qKZ map section} we define the qKZ equation and the qKZ map. In Section \ref{eqKZ section} we introduce the extended qKZ equation. Section \ref{main section} contains our main results, Theorem \ref{la1=k} and Corollaries \ref{qKZ same kind}-\ref{qKZ 2 to 1}. We describe the singularities of the hypergeometric map in Section \ref{singularities section}.

\section{Rational hypergeometric space of functions}\label{rational section}

\subsection{The Lie algebra $sl(2)$}\label{sl(2) section}

Let $e$, $f$, $h$ be generators of the Lie algebra $sl(2)$ such that 
\be
[h,e]=e, \qquad [h,f]=-f, \qquad [e,f]=2h.
\ee

For an $sl(2)$ module $M$, 
let $M^*$ be its restricted dual with an $sl(2)$ module structure defined by
\be
\langle e\varphi,x\rangle=\langle \varphi,fx\rangle,  \qquad
\langle f\varphi,x\rangle=\langle \varphi,ex\rangle,  \qquad
\langle h\varphi,x\rangle=\langle \varphi,hx\rangle,
\ee
for all $x\in M$, $\varphi \in M^*$. 

For a complex number $\la$, denote $V(\la)$ the $sl(2)$ Verma module with highest weight $\la$ and highest vector $v$. The vectors $\{f^lv,\; l\in\Z_{\geq 0}\}$ form a basis in $V(\la)$. The Shapovalov form on $V(\la)$ is a bilinear form $B_\la$ such that
\bean\label{Shapovalov}
B_\la(f^lv,f^lv)= l! \prod_{a=0}^{l-1}(2\la-a), \qquad 
B_\la(f^lv,f^kv)=0, \;l\neq k.
\eean
The Shapovalov form $B$ on a tensor product $V(\la_1)\T\dots\T V(\la_n)$ is defined by $B=B_{\la_1}\T\dots\T B_{\la_n}$. The dual map to the Shapovalov form defines a map of $sl(2)$ modules
\be
Sh:\; V(\la_1)\T\dots\T V(\la_n) \to (V(\la_1)\T\dots\T V(\la_n))^*\;.
\ee

For $\la\in\half\Z_{\geq 0}$, denote $L(\la)$ the ($2\la+1$)- dimensional irreducible $sl(2)$ module. 

\subsection{The Yangian $Y(gl(2))$}\label{Yangian}

The Yangian $Y(gl(2))$ is an associative algebra with an infinite set of generators $T_{i,j}^{(a)}$,  $i,j=1,2$, $a=0,1,\dots$, subject to the relations:
\be\label{yangian relation}
R(x-y)T_{(1)}(x)T_{(2)}(y)=T_{(2)}(y)T_{(1)}(x)R(x-y),
\ee
where $R(x)=(x\operatorname{Id}+P)\in\End(\C^2\T\C^2)$, $P\in\End(\C^2\T\C^2)$ is the operator of permutation of the two factors and $T_{ij}(u)=\Sum_{s=0}^\infty T_{ij}^{(a)}u^{-a}$ are generating series, 
$T_{(1)}(x)=1\T T(x)$, $T_{(2)}(x)=T(x)\T 1$.

In this paper we take tensor products and dual of Yangian modules using the comultiplication $\Dl\, :\,Y(gl(2))\to Y(gl(2))\T Y(gl(2))$ and the antipode $S\,:\,Y(gl(2))\to Y(gl(2))$ given by 
\bean\label{Dl}
\Dl\,:\, T_{ij}(u)\mapsto \sum_{a=1}^2 T_{ia}(u)\T T_{aj}(u),
\eean
\be
S\,:\,T_{ij}(u)\mapsto T_{ji}(u), \qquad i,j=1,2.
\ee 

For a complex number $z$, there is an automorphism $\rho_z\,:\, Y(gl(2))\to Y(gl(2))$ of the form
\be
\rho_z\,:\,T_{ij}(u)\mapsto T_{ij}(u-z),\qquad i,j=1,2.
\ee

The evaluation morphism $\epe\,:\, Y(gl(2))\to U(sl(2))$, to the universal enveloping algebra of $sl(2)$, has the form
\be
\epe\,:\,T_{11}(u)\mapsto h/u,\qquad \epe\,:\,T_{12}(u)\mapsto f/u,
\ee
\be
\epe\,:\,T_{21}(u)\mapsto e/u,\qquad \epe\,:\,T_{22}(u)\mapsto -h/u.
\ee

For complex numbers $z,\la$, denote $V(z,\la)$ the $sl(2)$ Verma module $V(\la)$ endowed with a structure of Yangian module via pull back with respect to the map $\epe\circ\rho_z$. The module $V(z,\la)$ is called  \emph{the evaluation Verma module}. 

For $\la\in\half\Z_{\geq 0}$, denote $L(z,\la)$ the $sl(2)$ module $L(\la)$ endowed with a structure of Yangian module via pull back with respect to the map $\epe\circ\rho_z$. The module $L(z,\la)$ is called  \emph{the evaluation finite dimensional module}. 

For more detail on the Yangian see \cite{CP}, \cite{TV}.

\subsection{The rational hypergeometric space}

Fix a natural number $n$ and $z=(\zn), \la=(\lan)\in\C^n$.

Define \emph{a rational hypergeometric space} $\F(z,\la)=\bigoplus\limits_{l=0}^\infty\,\F_l(z,\la)$, where $\F_l(z,\la)$ is the space of functions of the form
$$
P(\tell)\,\prod_{a=1}^n\,\prod_{b=1}^l\,\frac{1}{t_b-z_a-\la_a\!}\;
\prod_{1\leq a < b\leq l}\frac {t_a-t_b}{t_a-t_b+1}\;.
$$
Here $P$ is a polynomial with complex coefficients which is symmetric in variables  $\tell$ and has degree less than $n$ in each of variables $\tell$. 

\subsection{The Yangian action}\label{action section}

Let $f=f(\tell)$ be a function. For a permutation $\sigma\in{\Bbb S}^l$, define the functions
$[f]_{\sigma}^{rat}$  via the action of the simple transpositions 
$(i,i+1)\in {\Bbb S}^l$, $i=1,\dots,\l-1$, given by      

\be
[f]_{(i,i+1)}^{rat}(t_1,\dots, t_l)=
f(t_1,\dots,t_{i+1},t_i,\dots,t_l)\,\frac{t_i-t_{i+1}-1}{t_i-t_{i+1}+1}.
\ee

Let $T_{ij}(u)$, $i,j=1,2$, be the generating series for the Yangian $Y(gl(2))$
introduced in Section \ref{Yangian}. Set 
\be
\Ti_{ij}(u)\,=\,T_{ij}(u)\,\prod_{a=1}^n\,{u \over u-z_a-\la_a\!}\;,\qquad i,j=1,2\,.
\ee
Here the rational function in the right hand side is understood as its Laurent series expansion at $u=\infty$. In this paper we always use this convention in formulas of this kind. 

It is clear that the coefficients of the series $\Ti_{ij}(u)$ generate
$Y(gl(2))$. Following \cite{TV}, define an action of the coefficients of the series $\Ti_{ij}(u)$ in the hypergeometric space $\F(z,\la)$. Namely, for a function $f\in\F_l(z,\la)$, set:

\bean\label{yangian action}
(\Ti_{11}(u)f)(\tell)&=&f(\tell)\,
\prod_{a=1}^n\,\frac{u-z_a+\la_a}{ u-z_a-\la_a}\;
\prod_{a=1}^l\,\frac{u-t_a-1}{u-t_a}\ \,+  
\\ 
&+&\prod_{a=1}^l\,{u-t_a-1\over u-t_a}\;
\sum_{a=1}^l\,\left[\,{f(t_1,\dots,t_{l-1},u)\over u-t_l-1}\;
\prod_{b=1}^n\,{t_l-z_b+\la_b\over t_l-z_b-\la_b}\,\right]_{(a,l)}^{rat}\,,\notag
\eean

\bea
(\Ti_{22}(u)f)(\tell)&=&f(\tell)\,
\prod_{a=1}^l\,{u-t_a+1\over u-t_a}\,-
\\
&-&\,\prod_{a=1}^l\,{u-t_a+1\over u-t_a}\;\sum_{a=1}^l\,
\left[\,{f(u,t_2,\dots, t_l)\over u-t_1+1}\,\right]^{rat}_{(1,a)}\,,
\eea

\bea
(\Ti_{12}(u)f)(t_1,\dots,t_{l+1})&=&\sum_{a=1}^{l+1}\,\left[\,
{f(t_2,\dots, t_{l+1})\over u-t_1}\right.\times
\\ 
& \times & 
\left( \;\prod_{b=1}^n\,{t_1-z_b+\la_b\over t_1-z_b-\la_b}\;
\prod_{c=2}^{l+1}\>\>{u-t_c+1\over u-t_c}\,
{t_1-t_c-1\over t_1-t_c+1}\;-\right.
\\
&-&\left.\left.\,\prod_{b=1}^n\,{u-z_b+\la_b\over u-z_b-\la_b}\,
\prod_{c=2}^{l+1}\,{u-t_c-1\over u-t_c}\;\right)\>\right]^{rat}_{(1,a)}\,-
\eea
\be
-\,\prod_{a=1}^{l+1}\,{u-t_a+1\over u-t_a}
\sum_{\tsize{a,b=1\atop a\ne b}}^{l+1}\,
\left[\,{f(u,t_2,\dots,t_l)\over(u-t_1+1)\>(u-t_{l+1}+1)}\;
\prod_{c=1}^n\,{t_{l+1}-z_c+\la_c\over t_{l+1}-z_c-\la_c}\,\right]^{rat}_{(1,a)(b,l+1)}\,,
\ee
\be
(\Ti_{21}(u)f)(t_1,\dots, t_{l-1})\,=\,f(t_1,\dots, t_{l-1},u)\,
\prod_{a=1}^{l-1}\,{u-t_a-1\over u-t_a}\;,\qquad l>0\,,
\ee
and set $\Ti_{21}(u)f=0$ if $f\in\F_0(z,\la)$. Here $(1,a)$, $(a,l)$, $(b,l+1)$ are transpositions. 

By Lemma 4.2 in \cite{TV}, these formulas define a $Y(gl(2))$ module structure in the rational hypergeometric space $\F(z,\la)$ for any $z,\la\in\C^n$.

\subsection{Bases of the rational hypergeometric space}\label{bases section}

For natural numbers $n,l$, set $\Zb^n_l=\{\bar{l}=(l_1,\dots,\l_n)\in\Z^n_{\geq 0}\,|\, \sum\limits_{a=1}^nl_a=l\}$. For $\bar{l}\in\Zb^n_l$ 
and $i=0,1,\dots,n$, set $l^i=\sum\limits_{a=1}^il_a$.

For $\bar{l}\in\Zb_l^n$, define the \emph{rational weight function} 
$w_{\bar{l}}$ by 
\bea
w_{\bar{l}}(t,z,\la)=  
\sum_{\sigma \in {\Bbb S}^l}\left[\;\prod_{a=1}^n\frac{1}{l_a!}\prod_{b=l^{a-1}+1}^{l^a}
\left( \frac{1}{t_b-z_a-\la_a}\;\prod_{c=1}^{a-1}\,\frac{t_b-z_c+\la_c}{t_b-z_c-\la_c}
\right) \,\right] _\sigma^{rat}.
\eea

Denote  $V(z,\la)=V(z_1,\la_1)\T\dots\T V(z_n,\la_n)$ the tensor product of evaluation Verma modules.
The module $V(z,\la)$ has a basis given by monomials $\{f^{l_1}v_1\T\dots\T f^{l_n}v_n\}$, where $l_i\in\Z_{\geq 0}$ and $v_i$ are highest vectors of $V(z_i,\la_i)$. The dual space $V^*(z,\la)$ has the dual basis denoted by $(f^{l_1}v_1\T\dots\T f^{l_n}v_n)^*$.
Define a $\C$-linear map $\nu(z,\la)\,:\;V^*(z,\la)\;\to \;\F(z,\la)$
by the formula
\bean\label{nu}
\nu(z,\la)\,:\;(f^{l_1}v_1\T\dots\T f^{l_n}v_n)^*\mapsto w_{\bar{l}}(t,z,\la).
\eean

\begin{lemma}\label{isomor}(Cf. Lemma 4.5, Theorem 4.7 in \cite{TV}.) The map $\nu(z,\la)$ is a homomorphism of Yangian modules. Moreover, if $z_a-z_b+\la_a+\la_b\not\in\Z_{\geq 0}$ for all $a>b$, $a,b=1,\dots,n$, then the map $\nu(z,\la)$ is an isomorphism. $\;\Box$
\end{lemma} 

For $z,\la\in\C^n$ and permutations $\si,\si^\prime\in{\Bbb S}^n$, define $z_\si^{\si^\prime},\la_\si^{\si^\prime}\in\C^n$ by equations $(z_\si^{\si^\prime})_a+(\la_\si^{\si^\prime})_a=z_{\si(a)}+\la_{\si(a)}$ and
$(z_\si^{\si^\prime})_a-(\la_\si^{\si^\prime})_a=z_{\si^\prime(a)}-\la_{\si^\prime(a)}$, $a=1,\dots,n$.

\begin{lemma}\label{hyper spaces} For any permutations $\si,\si^\prime\in{\Bbb S}^n$, the identity map $\Id : \F(z_\si^{\si^\prime},\la_\si^{\si^\prime})\to \F(z,\la)$ is an isomorphism of Yangian modules.
\end{lemma}
\begin{proof}  The characteristic property of $z_\si^{\si^\prime},\la_\si^{\si^\prime}$ is 
\be
\prod_{a=1}^n(u-(z_\si^{\si^\prime})_a-(\la_\si^{\si^\prime})_a)=\prod_{a=1}^n(u-z_a-\la_a),\qquad \prod_{a=1}^n(u-(z_\si^{\si^\prime})_a+(\la_\si^{\si^\prime})_a)=\prod_{a=1}^n(u-z_a+\la_a).
\ee
The definition of the space $\F(z,\la)$ and of the Yangian action depends on $z,\la$ only through polynomials $ \prod\limits_{a=1}^n(u-z_a+\la_a)$ and  $\prod\limits_{a=1}^n(u-z_a-\la_a)$, see \Ref{yangian action}.
\end{proof}

\begin{lemma}\label{modules} Let $z_a-z_b+\la_a+\la_b\not\in\Z_{\geq 0}$, $a,b=1,\dots,n$. Then for any permutations $\si,\si^\prime\in{\Bbb S}^n$, the Yangian modules $V(z,\la)$ and $V(z_\si^{\si^\prime},\la_\si^{\si^\prime})$ are isomorphic.
\end{lemma}
{\it Proof 1\/}$:\;$
Under the assumption of the Lemma, the Shapovalov map $Sh:\,V(z_\si^{\si^\prime},\la_\si^{\si^\prime})\to V^*(z_\si^{\si^\prime},\la_\si^{\si^\prime})$ is an isomorpism of Yangian modules for any permutations $\si,\si^\prime\in{\Bbb S}^n$. Lemma \ref{modules} follows from Lemmas \ref{isomor} and \ref{hyper spaces}. 
{$\;\Box$ \par\vs}
{\it Proof 2\/}$:\;$ Both modules $V(z,\la)$ and $V(z_\si^{\si^\prime},\la_\si^{\si^\prime})$ are irreducible Yangian modules of highest weight
\be
\prod_{a=1}^n\;\frac{u-z_a+\la_a}{u-z_a-\la_a}\;.
\ee
For the definition of  highest weight and the classification of irreducible Yangian modules, see \cite{T}.
{$\;\Box$ \par\vs}

Let $z,\la$ be as in Lemma \ref{modules}.
Let $(\si^\prime\T\si)(z,\la): V(z,\la)\to V(z_\si^{\si^\prime},\la_\si^{\si^\prime})$ be the isomorphism of Yangian modules such that 
\be
(\si^\prime\T\si)(z,\la): v_1\T\dots\T v_n \mapsto v_1^\prime\T\dots\T v_n^\prime,
\ee
where $v_i, v_i^\prime$ are highest vectors generating $V(z_i,\la_i), V((z_\si^{\si^\prime})_i,(\la_\si^{\si^\prime})_i)$.

We have
\bean\label{si*si}
(\si^\prime\T\si)(z,\la)= (Sh(\la_\si^{\si^\prime}))^{-1}\circ (\nu(z_\si^{\si^\prime},\la_\si^{\si^\prime}))^{-1} \circ\nu(z,\la)\circ Sh(\la),
\eean
where the map $\nu$ is given by \Ref{nu}.

\begin{example}\label{R def} 
\rm
Let $n=2$.
Consider the transpositions $\si=\si^\prime=(1,2)$. Then 
\bean\label{R change}
z^{(1,2)}_{(1,2)}=(z_2,z_1),\qquad \la^{(1,2)}_{(1,2)}=(\la_2,\la_1).
\eean

The isomorpism $\widehat{R}(z,\la): V(z,\la)\to V(z^{(1,2)}_{(1,2)},\la^{(1,2)}_{(1,2)})$, given by formula \Ref{si*si}, will be called \emph{the rational $\widehat{R}$ matrix}. 
\end{example}

\begin{example}\label{N def} 
\rm
Let $n=2$.
Consider the case $\si=id$, $\si^\prime=(1,2)$. Then
\bean\label{N change}
(z^{(1,2)}_{\;id})_1=\half(z_1+z_2+\la_1-\la_2),\qquad
(z^{(1,2)}_{\;id})_2=\half(z_1+z_2-\la_1+\la_2),\notag
\\
(\la^{(1,2)}_{\;id})_1=\half(z_1-z_2+\la_1+\la_2),\qquad
(\la^{(1,2)}_{\;id})_2=\half(-z_1+z_2+\la_1+\la_2).
\eean
We call the isomorphism 
$N(z,\la): V(z,\la)\to V(z^{(1,2)}_{\;id},\la^{(1,2)}_{\;id})$, given by formula \Ref{si*si}, \emph{the rational $N$ matrix}. 
\end{example} 

\begin{example}\label{D def}
\rm
Let $n=2$.
Consider the case $\si^\prime=id$, $\si=(1,2)$. Then
\bean\label{D change}
(z_{(1,2)}^{\;id})_1=\half(z_1+z_2-\la_1+\la_2),\qquad
(z_{(1,2)}^{\;id})_2=\half(z_1+z_2+\la_1-\la_2),\notag
\\
(\la_{(1,2)}^{\;id})_1=\half(-z_1+z_2+\la_1+\la_2),\qquad
(\la_{(1,2)}^{\;id})_2=\half(z_1-z_2+\la_1+\la_2).
\eean
We call the isomorphism 
$D(z,\la): V(z,\la)\to V(z_{(1,2)}^{\;id},\la_{(1,2)}^{\;id})$, given by formula \Ref{si*si}, \emph{the rational $D$ matrix}. 
\end{example}

Our notations of $N$ and $D$ matrices come from the words ``numerator'' and ``denominator''.
We describe properties of the $\widehat{R}$, $D$ and $N$ matrices in the next Section.

\subsection{Properties of $R$, $D$ and $N$ matrices}

Let $n$ be a natural number. Let $z,\la\in\C^n$ be such that
$z_a-z_b+\la_a+\la_b\not\in\Z_{\geq 0}$, for all $a,b=1,\dots,n$.

\begin{lemma}\label{YBE} The group ${\Bbb S}^n\times {\Bbb S}^n$ transitively acts on the family of isomorphic Yangian modules $\{V(z_\si^{\si^\prime},\la_\si^{\si^\prime}), \si,\si^\prime\in{\Bbb S}^n\}$, $(\tau^\prime\times\tau)(z_\si^{\si^\prime},\la_\si^{\si^\prime}) :  V(z_\si^{\si^\prime},\la_\si^{\si^\prime})\to V(z_{\tau\si}^{\tau^\prime\si^\prime},\la_{\tau\si}^{\tau^\prime\si^\prime})$.
The action of generators is given by
\be
(a,a+1)\times id = N_{a,a+1},\qquad id\times (a,a+1) = D_{a,a+1},
\ee
where $N_{a,a+1}$, $D_{a,a+1}$ are $N$ and $D$ matrices acting in the $a$th and ($a+1$)st factors of the the corresponding tensor product.
We also have 
\be
N_{a,a+1}(z_{(a,a+1)}^{\;id},\la_{(a,a+1)}^{\;id})D_{a,a+1}(z,\la)=D_{a,a+1}(z^{(a,a+1)}_{\;id},\la^{(a,a+1)}_{\;id})N_{a,a+1}(z,\la)=\widehat{R}_{a,a+1}(z,\la).
\ee
\end{lemma}
\begin{proof} 
The Lemma follows from the definitions of $\widehat{R}$, $D$ and $N$ matrices. 
\end{proof}

Let $n=2$. Consider the $sl(2)$ Verma modules $V(\la_1), V(\la_2)$ of highest weights $\la_1,\la_2$ with highest vectors $v_1,v_2$. We have a decomposition of $sl(2)$ modules
\be
V(\la_1)\T V(\la_2)=\bigoplus_{l=0}^\infty V(\la_1+\la_2-l).
\ee
We choose singular vectors $v_{(\la_1,\la_2;\,l)}\in V(\la_1)\T V(\la_2)$ generating $V(\la_1+\la_2-l)$, so that in the standard basis the coefficient of $f^lv_1\T v_2$ is $(B_l(\la_1))^{-1}$:
\be
ev_{(\la_1,\la_2;\,l)}=0,\;\;\;
v_{(\la_1,\la_2;\,l)}=(B_l(\la_1))^{-1}\,f^lv_1\T v_2+\sum_{a=1}^l c^l_a(\la)f^{l-a}v_1\T f^av_2,\;\;\;l=0,1,\dots \,,
\ee
where $B_l(\la_1)=B_\la(f^lv,f^lv)$ is the value of the Shapovalov form, see \Ref{Shapovalov}, and $c^l_a(\la)$ are complex numbers.

\begin{thm}\label{spectral}(Spectral decomposition.) 
Let $A:V(z,\la)\to V(u,\om)$ be either $\widehat{R}$, $D$ or $N$ matrix, where $u,\om\in\C^2$ are the parameters of the target module. Then $A$ commutes with the $sl(2)$ action and maps
\be
A\,:\;v_{(\la_1,\la_2;\,l)}\;\mapsto \;\frac{B_l((\om^{(1,2)}_{\;id})_1)}{B_l((\la^{(1,2)}_{\;id})_1)}\;v_{(\om_1,\om_2;\,l)}, 
\ee
where $B_l(\la)=B_\la(f^lv,f^lv)$ is the value of the Shapovalov form, see  \Ref{Shapovalov}.
\end{thm}
\begin{proof}
In the case of $\widehat{R}$ matrix the spectral decomposition is well known, see Proposition 12.5.4 in \cite{CP} and formula (3.5) in \cite{TV}. In the two other cases the proof is similar.
\end{proof}

Let $n$ be any natural number, $z,\la\in \C^n$. We say that the parameters $z,\la\in\C^n$ are \emph{in a rational resonance} if the Yangian module $V(z,\la)$ is reducible. We say that the parameters $z,\la$ are \emph{in the first rational resonance} if there exists a unique pair of indices $a,b\in\{1,\dots,n\}$ such that $z_a-z_b+\la_a+\la_b=k\in\Z_{\geq 0}$ and if $a\leq b$. We say that the parameters $z,\la$ are \emph{in the second rational resonance} if there exists a unique pair of indices $a,b\in\{1,\dots,n\}$ such that $z_a-z_b+\la_a+\la_b=k\in\Z_{\geq 0}$ and if $a> b$.

Let $z,\la$ be in the first rational resonance. Then the Yangian module $V(z,\la)$ has a unique nontrivial submodule. The submodule is isomorphic to a tensor product $V(z^\prime,\la^\prime)$ of evaluation Verma modules with parameters $z^\prime,\la^\prime\in\C^n$ such that $(z^\prime,\la^\prime)\neq (z,\la)$. For example, if $2\la_1=k\in\Z_{\geq 0}$, i.e. $a=b=1$, then $z^\prime=z$, $\la^\prime=(-\la_1-1, \la_2,\dots,\la_n)$. 

Let $z,\la$ be in the second rational resonance. Then the Yangian module $V(z,\la)$ has a unique nontrivial submodule. The factor module is isomorphic to a tensor product $V(z^\prime,\la^\prime)$ of evaluation Verma modules with parameters $z^\prime,\la^\prime\in\C^n$ such that $(z^\prime,\la^\prime)\neq (z,\la)$. For example, if
$z_2-z_1+\la_1+\la_2=k\in\Z_{\geq 0}$, i.e. $a=2,\;b=1$, then $z^\prime=z^{(12)}_{\;id}$, $\la^\prime=(-(\la^{(12)}_{\;id})_1-1, (\la^{(12)}_{\;id})_2,\dots,(\la^{(12)}_{\;id})_n)$, see Section 12.1 in \cite{CP} and references therein.

\begin{lemma}\label{same kind}
Let $z,\la\in\C^n$ and $z_\si^{\si^\prime},\la_\si^{\si^\prime}\in\C^n$ be either both in the first rational resonance or both in the second rational resonance. Then the map 
$(\si^\prime\times\si)(z,\la): \; V(z,\la)\to V(z_\si^{\si^\prime},\la_\si^{\si^\prime})$ is a well defined isomorphism of Yangian modules. It maps the submodule $V(z^\prime,\la^\prime)\subset V(z,\la)$ onto the submodule $V((z_\si^{\si^\prime})^\prime,(\la_\si^{\si^\prime})^\prime)\subset V(z_\si^{\si^\prime},\la_\si^{\si^\prime})$, and the map $(\si^\prime\times\si)(z,\la)$ restricted to the submodule $V(z^\prime,\la^\prime)$ coincides with the map 
$(\si^\prime\times\si)(z^\prime,\la^\prime)$ up to a non-zero scalar multiplier, depending on the choice of the inclusions $V(z^\prime,\la^\prime)\hookrightarrow V(z,\la)$ and $V((z_\si^{\si^\prime})^\prime,(\la_\si^{\si^\prime})^\prime)\hookrightarrow V(z_\si^{\si^\prime},\la_\si^{\si^\prime})$.
\end{lemma}
\begin{proof}
Lemma \ref{same kind} follows from Theorem \ref{spectral}.
\end{proof}

\begin{lemma}\label{1 to 2}
Let $z,\la\in\C^n$ be in the first rational resonance
and $z_\si^{\si^\prime},\la_\si^{\si^\prime}\in\C^n$ be in the second rational resonance. Then the map 
$(\si^\prime\times\si)(z,\la): \; V(z,\la)\to V(z_\si^{\si^\prime},\la_\si^{\si^\prime})$ is a well defined homomorphism of Yangian modules. 
The kernel of this homomorphism is 
the submodule $V(z^\prime,\la^\prime)\subset V(z,\la)$ and the image is the proper submodule in $V(z_\si^{\si^\prime},\la_\si^{\si^\prime})$.
\end{lemma}
\begin{proof}
Lemma \ref{1 to 2} follows from Theorem \ref{spectral}.
\end{proof}

\begin{lemma}\label{2 to 1}
Let $z,\la\in\C^n$ be in the second rational resonance
and $z_\si^{\si^\prime},\la_\si^{\si^\prime}\in\C^n$ be in the first rational resonance. Then the map 
$(\si^\prime\times\si)(\tilde{z},\tilde{\la}): \; V(\tilde{z},\tilde{\la})\to V(\tilde{z}_\si^{\si^\prime},\tilde{\la}_\si^{\si^\prime})$ has a simple pole at $\tilde{z}=z,\tilde{\la}=\la$. The residue $Res:=\res_{\tilde{z}=z,\,\tilde{\la}=\la}(\si^\prime\times\si)(\tilde{z},\tilde{\la})$
is a well defined homomorphism of Yangian modules. 
The kernel of this homomorphism is 
the nontrivial submodule $U$ of $V(z,\la)$ and the image is the submodule 
$V((z_\si^{\si^\prime})^\prime,(\la_\si^{\si^\prime})^\prime)\subset V(z_\si^{\si^\prime},\la_\si^{\si^\prime})$. Thus, up to a scalar multiplier, depending on the choice of the factor map $V(z,\la)\to V(z^\prime,\la^\prime)\simeq V(z,\la)/U$, and the inclusion $V((z_\si^{\si^\prime})^\prime,(\la_\si^{\si^\prime})^\prime)\hookrightarrow V(z_\si^{\si^\prime},\la_\si^{\si^\prime})$, 
the map $Res$ defines a homomorphism $V(z^\prime,\la^\prime)\to V((z_\si^{\si^\prime})^\prime,(\la_\si^{\si^\prime})^\prime)$. The scalar multiplier can be chosen so that the map $Res$ coincides with the isomorphism $(\si^\prime\times\si)(z^\prime,\la^\prime)$.
\end{lemma}
\begin{proof}
Lemma \ref{2 to 1} follows from Theorem \ref{spectral}.
\end{proof}

\subsection{The first resonance in the rational hypergeometric space}\label{factor section} 

Let $z,\la\in\C^n$. Suppose $z,\la$ are in the first rational resonance. We have $z_a-z_b+\la_a+\la_b=k$ for some  $k\in\Z_{\geq 0}$, $a,b\in\{1,\dots,n\}$, $a\leq b$. Choose $\si,\si^\prime\in{\Bbb S}^n$, such that $\si(a)=\si^\prime(b)=1$. Then by Lemma \ref{same kind}, the map $\si^\prime\times\si$ is a well defined isomorphism and $2(\la_\si^{\si^\prime})_1=k$.

Let $z,\la\in\C^n$ and $2\la_1=k$, $k\in \Z_{\geq 0}$. Set $z^\prime=z$, $\la^\prime=(-\la_1-1,\la_2,\dots,\la_n)$.

Define a linear map $\iota^*(z,\la): \F(z,\la)\to \F(z^\prime,\la^\prime)$ as follows. For a function \linebreak
$f(\tell)\in\F_l(z,\la)$, set 
\be
\tilde{f}(\tell)=f(\tell)(t_1-z_1-\la_1)\;\prod_{a=1}^{k+1}\prod_{b=k+2}^l\;\frac{t_a-t_b+1}{t_a-t_b-1}, \qquad {\rm if}\; l>k,
\ee 
and $\tilde{f}=0$, if $l\leq k$.
Define 
\be
(\iota^* f)(t_1\dots,t_{l-k-1})= \tilde {f}(z_1+\la_1,z_1+\la_1-1,\dots,z_1-\la_1,t_1,\dots,t_{l-k-1}).
\ee
\begin{thm}\label{iota} The map $\iota^*(z,\la): \F(z,\la)\to \F(z^\prime,\la^\prime)$ is a surjective homomorphism of Yangian modules.
\end{thm}
\begin{proof}
The map $\iota^*$ is well defined, it is a surjection, see the definition of the rational hypergeometric space. 
 
If $t_a=z_1+\la_1-a+1$, $a=1,\dots, k+1$, then we have the following identities:
\be
\prod_{a=1}^{k+1}\,\frac{u-t_a+1}{u-t_a}\;\prod_{a=1}^n(u-z_a-\la_a)\,=\,
\prod_{a=1}^n \,(u-z^\prime_a-\la^\prime_a),
\ee
\be
\prod_{a=1}^{k+1}\,\frac{t_a-u+1}{t_a-u-1}\;\prod_{a=1}^n\,\frac{u-z_a+\la_a}{t-z_a-\la_a}\,=\,
\prod_{a=1}^n\,\frac{u-z^\prime_a+\la^\prime_a}{u-z^\prime_a-\la^\prime_a}.
\ee
It is a straightforward calculation to check that the map $\iota^*$ commutes with the Yangian action given by \Ref{yangian action}, using the above identities.
\end{proof}

Let $z,\la\in\C^n$ and $2\la_1=k\in \Z_{\geq 0}$. Assume 
$z_a-z_b+\la_a+\la_b\not\in\Z_{\geq 0}$ for all $a>b$, $a,b=1,\dots,n$. Then by Lemma \ref{isomor},
\be
\F(z,\la)\simeq V^*(z,\la)\simeq V^*(z_1,\la_1))\T V^*((z_2,\dots,z_n),(\la_2,\dots,\la_n)),
\ee
\be
\F(z^\prime,\la^\prime)\simeq V^*(z^\prime,\la^\prime)\simeq V^*(z^\prime_1,\la^\prime_1)\T V^*((z_2,\dots,z_n),(\la_2,\dots,\la_n)).
\ee 
The Yangian module $V^*(z_1,\la_1)$ has a submodule, isomorphic to ($k+1$)-dimensional irreducible evaluation module $L(z_1,\la_1)$. Then $\iota^*(z,\la)=\iota^*_1\T \Id$, where 
\be
\iota^*_1: V^*(z_1,\la_1)\to V^*(z_1,\la_1)/L(z_1,\la_1)\simeq V^*(z^\prime_1,\la^\prime_1)
\ee
is the factorization morphism mapping $f^{k+1}v_1\T v_2\T\dots\T v_n$ to $D(k)v^\prime_1\T\dots\T v_n^\prime$. The constant $D(k)$ is given by
\be
D(k)=\prod_{a=2}^{k+1} \frac{1}{t_a-z_1-\la_1}\,\prod_{1\leq a<b\leq k+1}\frac{t_a-t_b}{t_a-t_b+1},
\ee
where we set $t_a=z_1+\la_1-a+1$, $a=1,\dots, k+1$.

\section{Trigonometric hypergeometric space of functions}\label{trigonometric section}
\subsection{The quantum group $U_qsl(2)$}

Let $q$ be a complex number different from $\pm 1$. Let $p$ be a complex number such that $q=\exp(\pi i/p)$. We always assume $\Real p < 0$. For a complex number $x$, by $q^x$ we mean $\exp(\pi i x/p)$.

Let $e_q$, $f_q$, $q^h$, $q^{-h}$ be generators of $U_qsl(2)$ such that
\be
q^hq^{-h}=q^{-h}q^h=1, \qquad [e_q,f_q]=\frac{q^{2h}-q^{-2h}}{q-q^{-1}},
\qquad
q^he_q=qe_qq^h, \qquad  q^hf_q=q^{-1}f_qq^h.
\ee

A comultiplication $\Delta_q : U_qsl(2)\to U_qsl(2)\T U_qsl(2)$ is given by
\be
\Delta_q(q^{\pm h)}=q^{\pm h}\T q^{\pm h},
\qquad
\Delta_q(e_q)=e_q\T q^h + q^{-h}\T e_q,\qquad 
\Delta_q(f_q)=f_q\T q^h + q^{-h}\T f_q.
\ee
The comultiplication defines a module structure on tensor products of 
$U_qsl(2)$ modules.

An antipode $S_q:U_qsl(2)\to U_qsl(2)$ is given by
\be
S_q(e_q)=f_q,\qquad S_q(f_q)=e_q,\qquad S_q(q^{\pm h})=q^{\pm h}.
\ee
The antipode defines a module structure on a space dual to an $U_qsl(2)$ module.

For $\la \in \C$, denote $V_q(\la)$ the $U_qsl(2)$ Verma module with highest weight $q^\la$ and highest vector $v^q$. The vectors $\{f_q^lv^q,\; l\in\Z_{\geq 0}\}$ form a basis in $V_q(\la)$. The quantum Shapovalov form on $V_q(\la)$ is a bilinear form $B^q_\la$ such that
\bean\label{q Shapovalov}
B^q_\la\,(f_q^lv^q,f_q^lv^q)= [l]_q!\,[2\la-a]_q, \qquad 
B^q_\la\,(f_q^lv^q,f_q^kv^q)=0, \;l\neq k,
\eean
where $[n]_q=(q^n-q^{-n})/(q-q^{-1})=\sin (\pi n/p)/\sin(\pi/p)$ are the $q$-numbers, $[l]_q!=[1]_q![2]_q!\dots[l]_q!$.

The quantum Shapovalov form $B^q$ on a tensor product $V_q(\la_1)\T\dots\T V_q(\la_n)$ is defined by $B^q=B^q_{\la_1}\T\dots\T B^q_{\la_n}$. The dual map to the quantum Shapovalov form defines a map of $U_qsl(2)$ modules
\be
Sh_q:\; V_q(\la_1)\T\dots\T V_q(\la_n) \to (V_q(\la_1)\T\dots\T V_q(\la_n))^*\;.
\ee

For $\la\in\half\Z_{\geq 0}$, denote $L_q(\la)$ the $q$ deformation of $sl(2)$ module $L(\la)$. $L_q(\la)$ is a ($2\la+1$)-dimensional $U_qsl(2)$ module. 

\subsection{The Hopf algebra $\Ugh$}\label{Ugh}

The quantum affine algebra $\Ugh$ is a unital
associative algebra with generators $L_{ij}^{(+0)},\ L_{ji}^{(-0)}\!$,
$1\!\leq\!j\!\leq\!i\!\leq 2$, and $L_{ij}^{(s)}\!$, $i,j=1,2$, $s=\pm 1,\pm 2,\ldots$, subject to relations \Ref{Ughr}.

Let $e_{ij}\in\End(\C^2)$ $i,j=1,2$, be the standard matrix units with the only nonzero entry $1$ at the intersection of the $i$-th row and $j$-th column. Set

\bea
R(\xi) &=& (\xi q-q^{-1})\>(e_{11}\T e_{11}+e_{22}\T e_{22})\;+
\\
&+&(\xi-1)\>(e_{12}\T e_{12}+e_{21}\T e_{21})\,+\,
\xi(q-q^{-1})\>e_{12}\T e_{21}\,+\,(q-q^{-1})\>e_{21}\T e_{12}\,.
\eea

Introduce the generating series
$L^{\pm}_{ij}(u)=L_{ij}^{(\pm0)}+\sum\limits_{s=1}^\infty L_{ij}^{(\pm s)}u^{\pm s}\!$. The
relations in $\Ugh$ have the form

\bean\label{Ughr}
R(\xi/\zt)\>L^{\pm}_{(1)}(\xi)\>L^{\pm}_{(2)}(\zt)\,=\,
L^{\pm}_{(2)}(\zt)\>L^{\pm}_{(1)}(\xi)\>R(\xi/\zt)\,,
\\
R(\xi/\zt)\>L^+_{(1)}(\xi)\>L^-_{(2)}(\zt)\,=\,
L^-_{(2)}(\zt)\>L^+_{(1)}(\xi)\>R(\xi/\zt)\,,\notag
\eean
\be
L_{11}^{(\pm 0)}L_{22}^{(\pm 0)}\>=\>1\,,\qquad
L_{22}^{(\pm 0)}L_{11}^{(\pm 0)}\>=\>1\,,\qquad
L_{ii}^{(\pm 0)}L_{ii}^{(\mp 0)}\,=\,1\,,\qquad i=1,2\,,
\ee
where $L^\pm_{(1)}(\xi)=L^\pm(\xi)\T 1$ and
$L^\pm_{(2)}(\xi)=1\T L^\pm(\xi)$.

In this paper we take tensor products and dual of $\Ugh$ modules using the comultiplication $\widehat{\Dl}_q\, :\,\Ugh\to\Ugh\T\Ugh$ and the antipode $\widehat{S}_q\,:\,\Ugh\to \Ugh$ given by 
\bean\label{q Dl}
\widehat{\Dl}_q\,:\, L^\pm_{ij}(\xi)\mapsto \sum_{a=1}^2 L^\pm_{ai}(\xi)\T L^\pm_{ja}(\xi),
\eean
\be
\widehat{S}_q\,:\,L^\pm_{ij}(\xi)\mapsto L^\pm_{ji}(\xi),\qquad i,j=1,2.
\ee 
Note that our choice of comultiplication is, in a sense, opposite to the comultiplication we chose for Yangian, see \Ref{Dl}. 

For a complex number $\zt$, there is an automorphism $\rho^q_\zt\,:\, \Ugh\to \Ugh$ of the form
\be
\rho^q_\zt\,:\,L^\pm_{ij}(\zt)\mapsto L^\pm_{ij}(\xi/\zt),\qquad i,j=1,2.
\ee

The evaluation morphism $\epe^q\,:\, \Ugh\to U_qsl(2)$ has the form
\be
\epe^q\,:\; L^{\pm}_{11}(\xi)\,\mapsto\,q^{\mp h}-q^{\pm h}\xi^{\pm 1}\,,\qquad
\epe^q\,:\; L^{\pm}_{22}(\xi)\,\mapsto\,q^{\pm h}-q^{\mp h}\xi^{\pm 1}\,,
\ee
\be
\epe^q\,:\;L^+_{12}(\xi)\,\mapsto\,-f_q\,(q-q^{-1})\>\xi\,,\qquad
\epe^q\,:\;L^-_{12}(\xi)\,\mapsto\,f_q\,(q-q^{-1}),
\ee
\be
\epe^q\,:\;L^+_{21}(\xi)\,\mapsto\,-e_q\,(q-q^{-1})\,,\qquad
\epe^q\,:\;L^-_{21}(\xi)\,\mapsto\,e_q\,(q-q^{-1})\>\xi^{-1}\,.
\ee

For complex numbers $z,\la$, denote $V_q(z,\la)$ the $U_qsl(2)$ Verma module $V_q(\la)$ endowed with a structure of $\Ugh$ module via pull back with respect to the map $\epe^q\circ\rho^q_{q^{2z}}$. The module $V_q(z,\la)$ is called \emph{the quantum evaluation Verma module}. 

For $\la\in\half\Z_{\geq 0}$, denote $L_q(z,\la)$ the $U_qsl(2)$ module $L_q(\la)$ endowed with a structure of $\Ugh$ module via pull back with respect to the map $\epe^q\circ\rho^q_{q^{2z}}$. The module $L_q(z,\la)$ is called  \emph{the quantum evaluation finite dimensional module}. 

For more detail on the affine quantum group $\Ugh$ see \cite{CP}, \cite{TV}.

\subsection{The trigonometric hypergeometric space}

Fix a natural number $n$, $z=(\zn), \la=(\lan)\in\C^n$ and a complex number $q=\exp(\pi i /p)$.

Define \emph{a trigonometric hypergeometric space} $\F^q(z,\la)=\bigoplus\limits_{l=0}^\infty\,\F^q_l(z,\la)$, where $\F^q_l(z,\la)$ is the space of functions of the form
$$
P(q^{2t_1},\dots,q^{2t_n})\,\prod_{a=1}^n\,\prod_{b=1}^l\,\frac{q^{z_a-t_b}}{\sin(\pi(t_b-z_a-\la_a)/p)\!}\;
\prod_{1\leq a < b\leq l}\frac {\sin(\pi(t_a-t_b)/p)}{\sin(\pi(t_a-t_b+1)/p)}\;.
$$
Here $P$ is a polynomial with complex coefficients which is symmetric in variables  $q^{2t_1},\dots,q^{2t_n}$ and has degree less than $n$ in each of variables $q^{2t_1},\dots,q^{2t_n}$. 

\subsection{The $\Ugh$ action}

Let $f=f(\tell)$ be a function. For a permutation $\sigma\in{\Bbb S}^l$, define the functions
$[f]_{\sigma}^{trig}$  via the action of the simple transpositions 
$(i,i+1)\in {\Bbb S}^l$, $i=1,\dots,\l-1$, given by      

\be
[f]_{(i,i+1)}^{trig}(t_1,\dots, t_l)=
f(t_1,\dots,t_{i+1},t_i,\dots,t_l)\,\frac{\sin(\pi(t_i-t_{i+1}-1)/p)}{\sin(\pi(t_i-t_{i+1}+1)/p)}.
\ee

Let $L^\pm_{ij}(u)$, $i,j=1,2$, be the generating series for the quantum affine group $\Ugh$
introduced in Section \ref{Ugh}. Set 
\be
\Li^\pm_{ij}(\xi)\,=\,L^\pm_{ij}(\xi))\,\prod_{a=1}^n\,{\pm iq^{\pm(z_a-u)} \over 2\sin(\pi(u-z_a-\la_a)/p)\!}\;,\qquad i,j=1,2\,,
\ee
where $\xi=q^{2u}$. It is clear that the coefficients of the series $\Li^\pm_{ij}(u)$ generate
$\Ugh$. Following \cite{TV}, define an action of the coefficients of the series $\Li^\pm_{ij}(u)$ in the hypergemetric space $\F^q(z,\la)$. Namely, for a function $f\in\F^q_l(z,\la)$, set:

\bean\label{Ugh action}
\lefteqn{(\Li^\pm_{11}(\xi)f)(\tell)=}
\\ &&
=f(\tell)\,
\prod_{a=1}^n\,\frac{\sin(\pi(u-z_a+\la_a)/p)}{\sin(\pi( u-z_a-\la_a)/p)}\;
\prod_{a=1}^l\,\frac{\sin(\pi(u-t_a-1)/p)}{\sin(\pi(u-t_a)/p)}\ \,+ \sin(\pi/p)\times  \notag
\eean
\be
\times\prod_{a=1}^l\,{\sin(\pi(u-t_a-1)/p)\over \sin(\pi(u-t_a)/p)}\;
\sum_{a=1}^l\,\left[\,{f(t_1,\dots,t_{l-1},u)\;q^{u-t_l}\over \sin(\pi(u-t_l-1)/p)}\;
\prod_{b=1}^n\,{\sin(\pi(t_l-z_b+\la_b)/p)\over \sin(\pi(t_l-z_b-\la_b)/p)}\,\right]_{(a,l)}^{trig}\,,\notag
\ee

\bea
\lefteqn{
(\Li^\pm_{22}(\xi)f)(\tell)=f(\tell)\,
\prod_{a=1}^l\,{\sin(\pi(u-t_a+1)/p)\over \sin(\pi(u-t_a)/p)}\,-}
\\&&
-\,\sin(\pi/p)\prod_{a=1}^l\,{\sin(\pi(u-t_a+1)/p)\over \sin(\pi(u-t_a)/p)}\;\sum_{a=1}^l\,
\left[\,{f(u,t_2,\dots, t_l)q^{u-t_1}\over \sin(\pi(u-t_1+1)/p)}\,\right]^{trig}_{(1,a)}\,,
\eea

\bea
\lefteqn{
(\Li_{12}(\xi)f)(t_1,\dots,t_{l+1})\,=\,\sin(\pi/p)\,\sum_{a=1}^{l+1}\,\left[\,
{f(t_2,\dots, t_{l+1})\;q^{u-t_1}\over \sin(\pi(u-t_1)/p)}\right.\times}
\\ &&
\times\, 
\left( \;\prod_{b=1}^n\,{\sin(\pi(t_1-z_b+\la_b)/p)\over \sin(\pi(t_1-z_b-\la_b)/p)}\;
\prod_{c=2}^{l+1}\>\>{\sin(\pi(u-t_c+1)/p)\over \sin(\pi(u-t_c)/p)}\,
{\sin(\pi(t_1-t_c-1)/p)\over \sin(\pi(t_1-t_c+1)/p)}\;-\right.
\\&&
-\,\left.\left.\,\prod_{b=1}^n\,{\sin(\pi(u-z_b+\la_b)/p)\over \sin(\pi(u-z_b-\la_b)/p)}\,
\prod_{c=2}^{l+1}\,{\sin(\pi(u-t_c-1)/p)\over \sin(\pi(u-t_c)/p)}\;\right)\>\right]^{trig}_{(1,a)}\,-
\eea
\bea
&-&\sin^2(\pi/p)\prod_{a=1}^{l+1}\,{\sin(\pi(u-t_a+1)/p)\over \sin(\pi(u-t_a)/p)}
\left.\sum_{a,b=1,\; a\ne b}^{l+1}\,
\right[\,f(u,t_2,\dots,t_l)\times
\\
&\times&\left.
{q^{2u-t_1-t_{l+1}} \over \sin(\pi(u-t_1+1)/p)\>\sin(\pi(u-t_{l+1}+1)/p)}\;
\prod_{c=1}^n\,{\sin(\pi(t_{l+1}-z_c+\la_c)/p)\over \sin(\pi(t_{l+1}-z_c-\la_c)/p)}\,\right]^{trig}_{(1,a)(b,l+1)}\,,
\eea
\be
(\Li^\pm_{21}(\xi)f)(t_1,\dots, t_{l-1})\,=\,f(t_1,\dots, t_{l-1},u)\,
\prod_{a=1}^{l-1}\,{\sin(\pi (u-t_a-1)/p)\over \sin(\pi (u-t_a)/p)}\;,\qquad l>0\,,
\ee
and set $\Li^\pm_{21}(\xi)f=0$ if $f\in\F^q_0(z,\la)$. Here $(1,a)$, $(a,l)$, $(b,l+1)$ are transpositions, $\xi=q^{2u}$. 

By Lemma 4.15 in \cite{TV}, these formulas define a $\Ugh$ module structure in the trigonometric hypergeometric space $\F^q(z,\la)$ for any $z,\la\in\C^n$.

\subsection{Bases of the trigonometric hypergeometric space}

For $\bar{l}\in\Zb_l^n$, define the \emph{trigonometric weight function} 
$W_{\bar{l}}$ by 
\bea
\lefteqn{W_{\bar{l}}(t,z,\la)=}
\\
&&
\sum_{\sigma \in {\Bbb S}^l}
\left[\prod_{a=1}^n\prod_{d=1}^{l_m}
\frac{\sin(\pi/p)}{\sin(\pi d/p)}
\prod_{b=l^{a-1}+1}^{l^a}\!\!
\frac{\exp(\pi i(z_a-t_b)/p)}{\sin(\pi(t_b-z_a-\la_a)/p)}
\prod_{c=1}^{a-1}\frac{\sin(\pi(t_b-z_c+\la_c)/p)}{\sin(\pi(t_b-z_c-\la_c)/p)}
\right] _\sigma^{trig}.
\eea

Denote  $V_q(z,\la)=V_q(z_1,\la_1)\T\dots\T V_q(z_n,\la_n)$ the tensor product of evaluation Verma modules.
The module $V_q(z,\la)$ has a basis given by monomials $\{f_q^{l_1}v^q_1\T\dots\T f_q^{l_n}v^q_n\}$, where $l_i\in\Z_{\geq 0}$ and $v^q_i$ are highest vectors of $V_q(z_i,\la_i)$. The dual space $V_q^*(z,\la)$ has the dual basis denoted by $(f_q^{l_1}v^q_1\T\dots\T f_q^{l_n}v^q_n)^*$.
Define a $\C$-linear map $\nu_q(z,\la)\,:\;V_q^*(z,\la)\;\to \;\F^q(z,\la)$
by the formula
\bean\label{q nu}
\nu_q(z,\la)\,:\;(f_q^{l_1}v^q_1\T\dots\T f_q^{l_n}v^q_n)^*\mapsto \sin^l(\pi/p)\,W_{\bar{l}}(t,z,\la).
\eean

We often assume  $z_a-z_b+\la_a+\la_b\not\in \Z_{\geq 0}\oplus p\Z$, meaning that for any $k\in \Z_{\geq 0}, s\in\Z$, $z_a-z_b+\la_a+\la_b\neq  k+s$.

\begin{lemma}\label{q isomor}(Cf. Lemma 4.17, Theorem 4.19 in \cite{TV}.) The map $\nu_q(z,\la)$ is a homomorphism of $\Ugh$ modules. Moreover, if $z_a-z_b+\la_a+\la_b\not\in \Z_{\geq 0}\oplus p\Z$ for all $a>b$, $a,b=1,\dots,n$, then the map $\nu_q(z,\la)$ is an isomorphism. $\;\Box$
\end{lemma}

As before, for $z,\la\in\C^n$ and permutations $\si,\si^\prime\in{\Bbb S}^n$, define $z_\si^{\si^\prime},\la_\si^{\si^\prime}\in\C^n$ by equations
$(z_\si^{\si^\prime})_a+(\la_\si^{\si^\prime})_a=z_{\si(a)}+\la_{\si(a)}$ and
$(z_\si^{\si^\prime})_a-(\la_\si^{\si^\prime})_a=z_{\si^\prime(a)}-\la_{\si^\prime(a)}$, $a=1,\dots,n$.

\begin{lemma}\label{q hyper spaces} For any permutations $\si,\si^\prime\in{\Bbb S}^n$, the identity map $\Id : \F^q(z_\si^{\si^\prime},\la_\si^{\si^\prime})\to \F^q(z,\la)$ is an isomorphism of $\Ugh$ modules.
\end{lemma}
\begin{proof} (Cf. Lemma \ref{hyper spaces}.)
The definition of the space $\F^q(z,\la)$ and of the $\Ugh$ action depends on $z,\la$ only through quantities $ \prod\limits_{a=1}^n\sin(\pi(u-z_a+\la_a)/p)$ and  $\prod\limits_{a=1}^n\sin(\pi(u-z_a-\la_a)/p)$, see \Ref{Ugh action}.
\end{proof}

In what follows we assume that $q$ is not a root of unity, that is $p\not\in\Q$.

\begin{lemma}\label{q modules} Let $z_a-z_b+\la_a+\la_b\not\in \Z_{\geq 0}\oplus p\Z$, $a,b=1,\dots,n$. Assume $p\not\in\Q$. Then for any permutations $\si,\si^\prime\in{\Bbb S}^n$, the $\Ugh$ modules $V_q(z,\la)$ and $V_q(z_\si^{\si^\prime},\la_\si^{\si^\prime})$ are isomorphic.
\end{lemma}
\begin{proof} 
The proof of the Lemma is analogues to the proof of Lemma \ref{modules}.
\end{proof}

Let $z,\la,p$ be as in Lemma \ref{q modules}. Let $(\si^\prime\T\si)(z,\la): V_q(z,\la)\to V_q(z_\si^{\si^\prime},\la_\si^{\si^\prime})$ be the isomorphism of $\Ugh$ modules such that
\be
(\si^\prime\T\si)(z,\la): v^q_1\T\dots\T v^q_n \mapsto (v^q_1)^\prime\T\dots\T (v^q_n)^\prime,
\ee
where $v_i^q, (v_i^q)\prime$ are highest vectors generating $V_q(z_i,\la_i), V_q((z_\si^{\si^\prime})_i,(\la_\si^{\si^\prime})_i)$.

We have
\bean\label{q si*si}
(\si^\prime\T\si)(z,\la)= (Sh_q(\la_\si^{\si^\prime}))^{-1}\circ (\nu_q(z_\si^{\si^\prime},\la_\si^{\si^\prime}))^{-1} \circ\nu_q(z,\la)\circ Sh_q(\la),
\eean
where the map $\nu_q$ is given by \Ref{q nu}.

\begin{example}\label{q R def}
\rm
(Cf. Example \ref{R def}.) Let $n=2$.
Consider the transpositions $\si=\si^\prime=(1,2)$. Then 
$z^{\si ^\prime}_\si, \la^{\si ^\prime}_\si$ are given by \Ref{R change}.

The isomorpism $\widehat{R}^q(z,\la): V_q(z,\la)\to V_q(z^{(1,2)}_{(1,2)},\la^{(1,2)}_{(1,2)})$, given by formula \Ref{q si*si}, will be called \emph{the trigonometric $\widehat{R}^q$ matrix}. 
\end{example}

\begin{example}\label{q N def}
\rm
(Cf. Example \ref{N def}.)
Let $n=2$.
Consider the case $\si=id$, $\si^\prime=(1,2)$. Then
$z^{\si ^\prime}_\si, \la^{\si ^\prime}_\si$ are given by \Ref{N change}.

We call the isomorphism 
$N^q(z,\la): V_q(z,\la)\to V_q(z^{(1,2)}_{\;id},\la^{(1,2)}_{\;id})$, given by formula \Ref{q si*si}, \emph{the trigonometric $N^q$ matrix}. 
\end{example} 

\begin{example}\label{q D def}
\rm
(Cf. Example \ref{D def}.)
Let $n=2$.
Consider the case $\si^\prime=id$, $\si=(1,2)$. Then
$z^{\si ^\prime}_\si, \la^{\si ^\prime}_\si$ are given by \Ref{D change}.

We call the isomorphism 
$D^q(z,\la): V_q(z,\la)\to V_q(z_{(1,2)}^{\;id},\la_{(1,2)}^{\;id})$, given by formula \Ref{q si*si}, \emph{the trigonometric $D^q$ matrix}. 
\end{example}

We describe properties of the trigonometric $\widehat{R}^q$, $D^q$ and $N^q$ matrices in the next Section.

\subsection{Properties of $R^q$, $N^q$ and $D^q$ matrices}

Let $n$ be a natural number. Let $q=\exp(\pi i/p)$ be a complex number, not a root of unity. Let $z,\la\in\C^n$ be such that
$z_a-z_b+\la_a+\la_b\not\in \Z_{\geq 0}\oplus p\Z$, $a,b=1,\dots,n$.

\begin{lemma}\label{q YBE} The group ${\Bbb S}^n\times {\Bbb S}^n$ transitively acts on the family of isomorphic $\Ugh$ modules $\{V_q(z_\si^{\si^\prime},\la_\si^{\si^\prime}), \si,\si^\prime\in{\Bbb S}^n\}$, $(\tau^\prime\times\tau)(z_\si^{\si^\prime},\la_\si^{\si^\prime}) :  V_q(z_\si^{\si^\prime},\la_\si^{\si^\prime})\to V_q(z_{\tau\si}^{\tau^\prime\si^\prime},\la_{\tau\si}^{\tau^\prime\si^\prime})$.
The action of generators is given by
\be
(a,a+1)\times id = N^q_{a,a+1},\qquad id\times (a,a+1) = D^q_{a,a+1},
\ee
where $N^q_{a,a+1}$, $D^q_{a,a+1}$ are the trigonometric $N^q$ and $D^q$ matrices acting on the $a$th and ($a+1$)st factors.
We also have 
\be
N^q_{a,a+1}(z_{(a,a+1)}^{\;id},\la_{(a,a+1)}^{\;id})D^q_{a,a+1}(z,\la)=D^q_{a,a+1}(z^{(a,a+1)}_{\;id},\la^{(a,a+1)}_{\;id})N^q_{a,a+1}(z,\la)=\widehat{R}^q_{a,a+1}(z,\la).
\ee
\end{lemma}
\begin{proof} 
The Lemma follows from the definitions of $\widehat{R}^q$, $D^q$ and $N^q$ matrices. 
\end{proof}

Let $n=2$. Consider the $U_qsl(2)$ Verma modules $V_q(\la_1), V_q(\la_2)$ of highest weights $q^{\la_1}, q^{\la_2}$ with highest vectors $v^q_1,v^q_2$. We have a decomposition of $U_qsl(2)$ modules
\be
V_q(\la_1)\T V_q(\la_2)=\bigoplus_{l=0}^\infty V_q(\la_1+\la_2-l).
\ee
We choose singular vectors $v^q_{(\la_1,\la_2;\,l)}\in V_q(\la_1)\T V_q(\la_2)$ generating $V_q(\la_1+\la_2-l)$, so that in the standard basis the coefficient of $f_q^lv^q_1\T v^q_2$ is $(B_l^q(\la_1))^{-1}$:
\be
e_qv^q_{(\la_1,\la_2;\,l)}=0,\qquad v^q_{(\la_1,\la_2;\,l)}=(B_l^q(\la_1))^{-1}\,f_q^lv^q_1\T v^q_2+\sum_{a=1}^l d^l_a(\la)f_q^{l-a}v^q_1\T f_q^av^q_2,
\ee
$l=0,1,\dots$ ,
where $B_l^q(\la)=B^q_\la(f_q^lv^q,f_q^lv^q)$ is the value of the quantum Shapovalov form, see \Ref{q Shapovalov}, and $d^l_a(\la)$ are complex numbers.

\begin{thm}\label{q spectral}(Spectral decomposition.) 
Let $A^q:V_q(z,\la)\to V_q(u,\om)$ be either $\widehat{R}^q$, $D^q$ or $N^q$ matrix, where $u,\om\in\C^2$ are the parameters of the target module. Then $A^q$ commutes with the $U_qsl(2)$ action and maps
\be
A^q\,:\;v^q_{(\la_1,\la_2;\,l)}\;\mapsto \;\frac{B^q_l((\om_{(1,2)}^{\;id})_1)}{B^q_l((\la_{(1,2)}^{\;id})_1)}\;v^q_{(\om_1,\om_2;\,l)}, 
\ee
where $B^q_l(\la)=B^q_\la(f^lv,f^lv)$ is the value of the quantum Shapovalov form, see  \Ref{q Shapovalov}.
\end{thm}
\begin{proof}
In the case of $\widehat{R}^q$ matrix the spectral decomposition is well known, see Proposition 12.5.6 in \cite{CP} and formula (3.16) in \cite{TV}. In the two other cases the proof is similar.
\end{proof}

Note the difference in Theorems \ref{spectral} and \ref{q spectral} due to the opposite choice of the comultiplications, see \Ref{Dl} and \Ref{q Dl}.

Let $n$ be any natural number, $z,\la\in \C^n$. Let $q=\exp(\pi i/p)\in\C$ be not a root of unity. We say that the parameters $z,\la\in\C^n$ are \emph{in a trigonometric resonance} if the $\Ugh$ module $U_q(z,\la)$ is reducible. We say that the parameters $z,\la\in\C^n$ are \emph{in the first trigonometric resonance} if there exists a unique pair of indices $a,b\in\{1,\dots,n\}$ such that $z_a-z_b+\la_a+\la_b\in\Z_{\geq 0}\oplus p\Z$
and $a\geq b$. We say that the parameters $z,\la$ are \emph{in the second trigonometric resonance} if there exists a unique pair of indices $a,b\in\{1,\dots,n\}$ such that $z_a-z_b+\la_a+\la_b\in\Z_{\geq 0}\oplus p\Z$ and $a<b$.

Let $z,\la$ be in the first trigonometric resonance. Then the $\Ugh$ module $V_q(z,\la)$ has a unique nontrivial submodule. The submodule is isomorphic to a tensor product $V_q(z^\prime,\la^\prime)$ of quantum evaluation Verma modules with parameters $z^\prime,\la^\prime\in\C^n$, $(z^\prime,\la^\prime)\neq (z,\la)$. For example, if $2\la_1=k\in\Z_{\geq 0}$, i.e. $a=b=1$, then $z^\prime=z$, $\la^\prime=(-\la_1-1, \la_2,\dots,\la_n)$. 

Let $z,\la$ be in the second trigonometric resonance. Then the $\Ugh$ module $V_q(z,\la)$ has a unique nontrivial submodule. The factor module is isomorphic to a tensor product $V_q(z^\prime,\la^\prime)$ of quantum evaluation Verma modules with parameters $z^\prime,\la^\prime\in\C^n$, $(z^\prime,\la^\prime)\neq (z,\la)$. For example, if
$z_1-z_2+\la_1+\la_2=k\in\Z_{\geq 0}$, i.e. $a=1,\;b=2$, then $z^\prime=z_{(12)}^{\;id}$, $\la^\prime=(-(\la_{(12)}^{\;id})_1-1, (\la_{(12)}^{\;id})_2,\dots,(\la_{(12)}^{\;id})_n)$, see Section 12.2 in \cite{CP} and references therein.

\begin{lemma}\label{q same kind}
Let $z,\la\in\C^n$ and $z_\si^{\si^\prime},\la_\si^{\si^\prime}\in\C^n$ be either both in the first trigonometric resonance or both in the second trigonometric resonance. Then the map 
$(\si^\prime\times\si)(z,\la): \; V_q(z,\la)\to V_q(z_\si^{\si^\prime},\la_\si^{\si^\prime})$ is a well defined isomorphism of $\Ugh$ modules. In particular, it maps the submodule $V_q(z^\prime,\la^\prime)\subset V_q(z,\la)$ onto the submodule $V_q((z_\si^{\si^\prime})^\prime,(\la_\si^{\si^\prime})^\prime)\subset V_q(z_\si^{\si^\prime},\la_\si^{\si^\prime})$. Moreover, the map $(\si^\prime\times\si)(z,\la)$ restricted to the submodule $V_q(z^\prime,\la^\prime)$ coincides with the map 
$(\si^\prime\times\si)(z^\prime,\la^\prime)$ up to a non-zero scalar multiplier depending on the choice of the inclusions $V_q(z^\prime,\la^\prime)\hookrightarrow V_q(z,\la)$ and $V_q((z_\si^{\si^\prime})^\prime,(\la_\si^{\si^\prime})^\prime)\hookrightarrow V_q(z_\si^{\si^\prime},\la_\si^{\si^\prime})$.
\end{lemma}
\begin{proof}
Lemma \ref{q same kind} follows from Theorem \ref{q spectral}.
\end{proof}

\begin{lemma}\label{q 1 to 2}
Let $z,\la\in\C^n$ be in the first trigonometric resonance
and $z_\si^{\si^\prime},\la_\si^{\si^\prime}\in\C^n$ be in the second trigonometric resonance. Then the map 
$(\si^\prime\times\si)(z,\la): \; V(z,\la)\to V(z_\si^{\si^\prime},\la_\si^{\si^\prime})$ is a well defined homomorphism of $\Ugh$ modules. 
The kernel of this homomorphism is 
the submodule $V_q(z^\prime,\la^\prime)\subset V(z,\la)$ and the image is the proper submodule in $V_q(z_\si^{\si^\prime},\la_\si^{\si^\prime})$.
\end{lemma}
\begin{proof}
Lemma \ref{q 1 to 2} follows from Theorem \ref{q spectral}.
\end{proof}

\begin{lemma}\label{q 2 to 1}
Let $z,\la\in\C^n$ be in the second trigonometric resonance
and $z_\si^{\si^\prime},\la_\si^{\si^\prime}\in\C^n$ be in the first trigonometric resonance. Then the map 
$(\si^\prime\times\si)(\tilde{z},\tilde{\la}): \; V_q(\tilde{z},\tilde{\la})\to V_q(\tilde{z}_\si^{\si^\prime},\tilde{\la}_\si^{\si^\prime})$ has a simple pole at $\tilde{z}=z,\tilde{\la}=\la$. The residue $Res:=\res_{\tilde{z}=z,\,\tilde{\la}=\la}(\si^\prime\times\si)(\tilde{z},\tilde{\la})$
is a well defined homomorphism of $\Ugh$ modules. 
The kernel of this homomorphism is 
the nontrivial submodule $U_q$ of $V_q(z,\la)$ and the image is the submodule 
$V_q((z_\si^{\si^\prime})^\prime,(\la_\si^{\si^\prime})^\prime)\subset V_q(z_\si^{\si^\prime},\la_\si^{\si^\prime})$. Thus, up to a scalar multiplier, depending on the choice of the factor map $V_q(z,\la)\to V_q(z^\prime,\la^\prime)\simeq V_q(z,\la)/U_q$, and the inclusion $V_q((z_\si^{\si^\prime})^\prime,(\la_\si^{\si^\prime})^\prime)\hookrightarrow V_q(z_\si^{\si^\prime},\la_\si^{\si^\prime})$, the map $Res$ defines a homomorphism $V_q(z,\la)/U_q\simeq V_q(z^\prime,\la^\prime)\to V_q((z_\si^{\si^\prime})^\prime,(\la_\si^{\si^\prime})^\prime)$. The scalar multiplier can be chosen so that the map $Res$ coincides with the isomorphism 
$(\si^\prime\times\si)(z^\prime,\la^\prime)$.
\end{lemma}
\begin{proof}
Lemma \ref{q 2 to 1} follows from Theorem \ref{q spectral}.
\end{proof}

\subsection{The first resonance in the trigonometric hypergeometric space}
Let $z,\la\in\C^n$. Let $q=\exp(\pi i /p)$ be a complex number, not a root of unity. Suppose $z,\la$ are in the first trigonometric resonance.
We have  $z_a-z_b+\la_a+\la_b=k$ for some  $k\in\Z_{\geq 0}$, $a,b\in\{1,\dots,n\}$, $a\geq b$. Choose $\si,\si^\prime\in{\Bbb S}^n$, such that $\si(a)=\si^\prime(b)=1$. Then, by Lemma \ref{q same kind}, the map $\si^\prime\times\si$ is a well defined isomorphism and $2(\la_\si^{\si^\prime})_1=k$.

Let $z,\la\in\C^n$ and $2\la_1=k$, $k\in \Z_{\geq 0}$. As before, set $z^\prime=z$, $\la^\prime=(-\la_1-1,\la_2,\dots,\la_n)$.

Define a linear map $\iota_q^*(z,\la): \F^q(z,\la)\to \F^q(z^\prime,\la^\prime)$ as follows. For a function \linebreak
$f_q(\tell)\in\F^q_l(z,\la)$, set 
\be
\tilde{f}_q(\tell)=f_q(\tell)\,\sin(\pi(t_1-z_1-\la_1)/p)\;\prod_{a=1}^{k+1}\prod_{b=k+2}^l\;\frac{\sin(\pi(t_a-t_b+1)/p)}{\sin(\pi(t_a-t_b-1)/p)}, \; l>k,
\ee 
and $\tilde{f}=0$ if $l\leq k$. Define 
\be
(\iota_q^* f_q)(t_1\dots,t_{l-k-1})= \tilde {f}_q(z_1+\la_1,z_1+\la_1-1,\dots,z_1-\la_1,t_1,\dots,t_{l-k-1}).
\ee
\begin{thm}\label{q iota} The map $\iota_q^*(z,\la): \F^q(z,\la)\to \F^q(z^\prime,\la^\prime)$ is a surjective homomorphism of $\Ugh$ modules.
\end{thm}
\begin{proof}
(Cf. Theorem \ref{iota}.)
The map $\iota_q^*$ is well defined, it is a surjection, see the definition of the trigonometric hypergeometric space. 
 
If $t_a=z_1+\la_1-a+1$, $a=1,\dots, k+1$, then we have the following identities:
\be
\prod_{a=1}^{k+1}\,\frac{\sin(\pi(u-t_a+1)/p)}{\sin(\pi(u-t_a)/p)}\;\prod_{a=1}^n\,\sin(\pi(u-z_a-\la_a)/p)\,=\,
\prod_{a=1}^n\, \sin(\pi(u-z^\prime_a-\la^\prime_a)/p),
\ee
\be
\prod_{a=1}^{k+1}\,\frac{\sin(\pi(t_a-u+1)/p)}{\sin(\pi(t_a-u-1)/p)}\;\prod_{a=1}^n\,\frac{\sin(\pi(u-z_a+\la_a)/p)}{\sin(\pi(t-z_a-\la_a)/p)}\,=\,
\prod_{a=1}^n \,\frac{\sin(\pi(u-z^\prime_a+\la^\prime_a)/p)}{\sin(\pi(u-z^\prime_a-\la^\prime_a)/p)}.
\ee
It is a straightforward calculation to check that the map $\iota_q^*$ commutes with the $\Ugh$ action given by \Ref{Ugh action}, using the above identities.
\end{proof}

Let $z,\la\in\C^n$ and $2\la_1=k\in \Z_{\geq 0}$. Assume 
$z_a-z_b+\la_a+\la_b\not\in\Z_{\geq 0}\oplus p\Z$ for all $a>b$, $a,b=1,\dots,n$. Then by Lemma \ref{q isomor},
\be
\F^q(z,\la)\simeq V_q^*(z,\la)\simeq V_q^*(z_1,\la_1))\T V_q^*((z_2,\dots,z_n),(\la_2,\dots,\la_n)),
\ee
\be
\F^q(z^\prime,\la^\prime)\simeq V_q^*(z^\prime,\la^\prime)\simeq V_q^*(z^\prime_1,\la^\prime_1)\T V_q^*((z_2,\dots,z_n),(\la_2,\dots,\la_n)).
\ee 
The $\Ugh$ module $V_q^*(z_1,\la_1)$ has a submodule, isomorphic to ($k+1$)-dimensional irreducible evaluation module $L_q(z_1,\la_1)$. Then $\iota_q^*(z,\la)=(\iota_q^*)_1\T \Id$, where 
\be
(\iota_q^*)_1: V_q^*(z_1,\la_1)\to V_q^*(z_1,\la_1)/L_q(z_1,\la_1)\simeq V_q^*(z^\prime_1,\la^\prime_1)
\ee
is the factorization morphism, mapping $f_q^{k+1}v_1^q\T v_2^q\T\dots\T v_n^q$ to $D_q(k)(v^q_1)^\prime\T\dots (v^q_n)^\prime$. The constant $D_q(k)$ is given by
\be
D_q(k)=\sin^{k+1}(\pi/p)\prod_{a=2}^{k+1} \frac{1}{\sin(\pi(t_a-z_1-\la_1)/p)}\,\prod_{1\leq a<b\leq k+1}\frac{\sin(\pi(t_a-t_b)/p)}{\sin(\pi(t_a-t_b+1)/p)},
\ee
where we set $t_a=z_1+\la_1-a+1$, $a=1,\dots, k+1$.

\section{Hypergeometric pairing and its properties}

\subsection{The qKZ equation}\label{qKZ section}
Fix complex numbers $q=\exp(\pi i/p),\kappa=\exp(\mu)$, where we assume $\Real p <0$, $0<\Imag\mu<2\pi$. Let $\la\in\C^n$.
The rational $\widehat{R}$ matrix
$\widehat{R}((z_i,z_j),(\la_i,\la_j))$, depends on $z_i,z_j$ only through the difference $z_i-z_j$, see Theorem \ref{spectral}. The operator 
$R(z_i-z_j)=P\widehat{R}((z_i,z_j),(\la_i,\la_j))\in\End(V(\la_i)\T V(\la_j))$, where $P$ is the operator of permutation of the factors, will be called the rational $R$ matrix. 

The rational quantized Knizhnik-Zamolodchikov equation (qKZ) with values in a tensor product of $sl(2)$ Verma modules $V(\la_1)\T\dots\T V(\la_n)$
is a system of difference equations for a function $\Psi(z_1,...,z_n)\in V(\la_1)\T\dots\T V(\la_n)$.
The system of equations has the form
\be
\Psi(z_1,\dots,z_m+p,\dots,z_n)= H_m(z,\la)\;\Psi(z_1,\dots,z_n), \qquad m=1,\dots,n,
\ee
\bean\label{qKZ operators}
H_m(z,\la) &=&
R_{m, {m-1}}(z_m-z_{m-1}+p)\dots
R_{m,1}(z_m-z_1+p)\times\notag
\\ 
& \times &
\kappa^{-h_m} R_{m,n}(z_m-z_n)\dots R_{m,{m+1}}(z_m-z_{m+1}),
\eean
where $h_m$ is the operator $h\in sl(2)$ acting in the $m$-th factor,
$R_{i,j}(z_i-z_j)$ is the rational $R$ matrix acting in the $i$-th and $j$-th factors of the tensor product.

\subsection{The hypergheometric pairing}\label{qKZ map section}
Let $z=(z_1,\dots,z_n)\in\C^n,\;\la=(\la_1,\dots,\la_n)\in\C^n ,\; 
t=(t_1,\dots,t_l)\in\C^l$. Let $l\in\Z_{\geq 0}$.
The \emph{phase function} is defined by the following formula:
\be
\Phi_l(t,z,\la)=
\frac{1}{l!}\,
\exp(\mu\sum_{a=1}^lt_a/p)\prod_{a=1}^n\prod_{b=1}^l
\frac{\Gamma((t_b-z_a+\la_a)/p)}{\Gamma((t_b-z_a-\la_a)/p)}
\prod_{1\leq a< b\leq l}\frac{\Gamma((t_a-t_b-1)/p)}{\Gamma((t_a-t_b+1)/p)}.
\ee

Assume that the parameters $z,\la\in\C^n$ satisfy the condition $\Real (z_i+\la_i)<0$ and
$\Real (z_i-\la_i)>0$ for all $i=1,\dots,n$. For functions 
$w(t,z,\la)\in\F_l(z,\la)$,  and  
$W(t,z,\la)\in\F_l^q(z,\la)$, define the \emph{hypergeometric integral} $I(w,W)(z,\la)$ by the formula
\bean\label{int specific}
I(w,W)(z,\la)=\int\limits_{\Real t_i=0,\atop
i=1,\dots,l}\Phi_l(t,z,\la)w(t,z,\la)W(t,z,\la)\,d^lt,
\eean
where $d^lt=dt_1\ldots dt_l$.

The hypergeometric integral for generic $z,\la$ is defined by analytic continuation with respect to $z,\la$ and has the form
\be\label{int generic}
I(w,W)(z,\la)=\int\limits_{{\mathcal K}(z,\la;p)}\Phi_l(t,z,\la)w(t,z,\la)W(t,z,\la)\,d^lt,
\ee
for some suitable contour of integration ${\mathcal K}(z,\la;p)\subset\C^l$, see \cite{MV}. The hypergeometric integral $I(w,W)(z,\la)$ is a univalued meromorphic function of variables $z,\la$.

Define \emph{the hypergeometric pairing} ${\mathcal H}(z,\la): \F^q(z,\la)\T\F(z,\la)\,\to\,\C$ by
\be
{\mathcal H}(z,\la)\,:\;W(t,z,\la)\T w(t,z,\la)\mapsto I(w,W)(z,\la), 
\ee
if $w(t,z,\la)\in\F_l(z,\la)$, $W(t,z,\la)\in\F_l^q(z,\la)$, and
let  ${\mathcal H}(W\T w)=0$ if 
$w\in\F_l(z,\la)$, $W\in\F_k^q(z,\la)$, $k\neq l$.

Assume $p\not\in\Q$ and $z_a-z_b+\la_a+\la_b\not\in\Z_{\geq 0}+p\Z$, $a,b=1,\dots n$. Then by Lemmas \ref{isomor}, \ref{q isomor}, we have a map
\be
{\mathcal H}(z,\la)\,:\; V^*_q(z,\la)\T V^*(z,\la)\to \C.
\ee
Define \emph{the qKZ map}, 
\bean
{\rm qKZ}(z,\la) &:& V_q(z,\la)\to V(z,\la)\notag,
\\
{\rm qKZ}(z,\la) &:&
f_q^{l_1}v_1^q\T\dots\T f_q^{l_n}v_n^q\mapsto 
\Psi_{\bar{l}}(z,\la),\notag\label{qKZ map}
\\
\Psi_{\bar{l}}(z,\la)
&=&\sum_{\bar{m}\in\Zb_n^l} B^q_{\bar{l}}(\la)\,I(w_{\bar{m}},W_{\bar{l}})(z,\la)
\;f^{m_1}v_1\T\dots\T f^{m_n}v_n,\label{Psi}
\eean
where $B^q_{\bar{l}}(\la)=B^q(f_q^{l_1}v^q_1\T\dots\T f_q^{l_n}v_n^q,f_q^{l_1}v_1^q\T\dots\T f_q^{l_n}v_n^q)$ is the value of the quantum Shapovalov form, see \Ref{q Shapovalov}.
 
We have ${\rm qKZ}(z,\la)=\tilde{\mathcal H}(z,\la)\circ Sh_q(\la)$, where $\tilde{\mathcal H}(z,\la):V^*_q(z,\la)\to V(z,\la)$ is the map dual to the hypergeometric pairing.

The meromorphic functions  $\{\Psi_{\bar{l}}(z,\la), \; \bar{l}\in\Z^n_{\geq 0}\}$ form a  
basis of solutions of the rational qKZ equation, see Corollary 5.25 in \cite{TV}. In particular, if
$p\not\in\Q$ and $z_a-z_b+\la_a+\la_b\not\in\Z_{\geq 0}+p\Z$, $a,b=1,\dots n$, then
 the qKZ map is an isomorphism of vector spaces, meromorphic with respect to variables $z,\la$.

\subsection{The qKZ map and the $\widehat{R}$, $D$, $N$ matrices}\label{eqKZ section}

%
%

\begin{thm}\label{extended monodromy} 
Let $\Psi_{\bar{l}}(z,\la)$ be a hypergeometric solution of the qKZ equation with values in $V(z,\la)$ given by \Ref{Psi}. Then for any permutations $\si,\si^\prime\in{\Bbb S}^n$, the function $(\si^\prime\times\si)\, (\Psi_{\bar{l}}(z,\la))$ is a solution of the qKZ equation with values in $V (z^{\si^\prime}_\si,\la^{\si^\prime}_\si)$. Moreover, the following diagram is commutative:
\bean 
\renewcommand{\arraystretch}{1.5}
\begin{array}{ccc}
{\qquad V_q(z,\la)\qquad}&
\stackrel{{\rm qKZ}(z,\la)}
{\longrightarrow}
&{\qquad V(z,\la)\qquad}
\\
{\qquad\downarrow \;\si^\prime\times\si}&&{\qquad\downarrow\;\si^\prime\times\si}
\\
{V_q(z_\si^{\si^\prime},\la_\si^{\si^\prime})}&
\stackrel{{\rm qKZ}(z_\si^{\si^\prime},\la_\si^{\si^\prime})}
{\longrightarrow}
&{V(z_\si^{\si^\prime},\la_\si^{\si^\prime}).}
\end{array}
\eean
\end{thm}
\begin{proof}
We have $\Phi_l(t,z,\la)=\Phi_l(t,z_\si^{\si^\prime},\la_\si^{\si^\prime})$.
Theorem \ref{extended monodromy} follows from the definitions of all participating maps and Lemmas \ref{hyper spaces} and \ref{q hyper spaces}.
\end{proof}

\subsection{The extended qKZ equation}

Consider the space $\C^{2n}$ with coordinates $(z,\la)$, $z,\la\in\C^n$.
For permutations $\si^\prime,\si\in{\Bbb S}^n$, introduce shifts $T_\si^{\si^\prime}(p):\, \C^{2n}\to\C^{2n}$ by the formula
\be
T_\si^{\si^\prime}(p)(z,\la))=((\si^\prime)^{-1}\times (\si)^{-1})\circ T(p)\circ (\si^\prime\times\si) \,(z,\la),\qquad z,\la\in\C^n,
\ee
where $(\si^\prime\times\si) \,(z,\la)=(z_\si^{\si^\prime},\la_\si^{\si^\prime})$, and $T(p)(z,\la)=((z_1+p,z_2,\dots,z_l),(\la_1,\dots,\la_n))$.

Introduce a system of difference equations for a function $\Psi(z,\la)\in V(z,\la)$ by
\bean\label{eqKZ equation}
\Psi(T_\si^{\si^\prime}(p)(z,\la))=H_\si^{\si^\prime}(z,\la;p,\mu)\Psi(z,\la),
\eean
where
\be\label{eqKZ operators}
H_\si^{\si^\prime}(z,\la;p,\mu)= ((\si^\prime)^{-1}\times (\si)^{-1})\circ H_1(z_\si^{\si^\prime},\la_\si^{\si^\prime};p,\mu)\circ (\si^\prime\times\si)\in\End (V(z,\la)),\notag
\ee
and $H_1(z,\la;p,\mu)$ is the first qKZ operator, see \Ref{qKZ operators}. We call the system \ref{eqKZ equation} \emph{the extended qKZ equation}.

\begin{example}\label{qKZ in eqKZ}
\rm
The system \Ref{eqKZ equation} contains the qKZ equation \Ref{qKZ operators}.
Namely, for permutations $\si, \si^\prime$ such that  $\si=\si^\prime$, $\si(1)=i$ the equation \Ref{eqKZ equation} takes the form
$\Psi(z_1,\dots,z_i+p,\dots,z_n,\la)=H_i(z,\la)\Psi(z,\la)$, where $H_i(z,\la)$ is given by \Ref{qKZ operators}.
\end{example}

\begin{example}
\rm
The shifts $T_\si^{\si^\prime}(p): \C^{2n}\to\C^{2n}$ preserve the quantity $\la_1+\dots+\la_n$ and generate a group acting in the space $\C^{2n}$ of parameters $(z,\la)$ isomorpfic to $\Z^{2n-1}$. For example, we have
\be
T^{\;id}_{(1,2)}(z,\la)=((z_1+p/2,\,z_2+p/2,\,z_3,\dots,z_n),(\la_1-p/2,\,\la_2+p/2,\,\la_3,\dots,\la_n)). 
\ee
\end{example}

\begin{corollary}\label{extended qKZ}
The extended qKZ equation is a compatible (holonomic) system of difference equations, i.e. for any permutations $\si,\si^\prime,\tau,\tau^\prime\in{\Bbb S}^n$,
\be
H^{\si^\prime}_\si(T_\tau^{\tau^\prime}(p)(z,\la))H_\tau^{\tau^\prime}(z,\la)=
H_\tau^{\tau^\prime}(T^{\si^\prime}_\si(p)(z,\la))H^{\si^\prime}_\si(z,\la).
\ee
Moreover, the meromorphic functions $\{\Psi_{\bar{l}}(z,\la), \;\bar{l}\in\Z^n_{\geq 0}\}$, where the function $\Psi_{\bar{l}}(z,\la)$ is given by \Ref{Psi}, form a basis of solutions of the extended qKZ equation. 
\end{corollary}
\begin{proof}
Corollary \ref{extended qKZ} follows from Theorem \ref{extended monodromy}. 
\end{proof}

It would be interesting to find an algebraic proof of Corollary \ref{extended qKZ}.

\subsection{The qKZ map, factormodules and submodules}\label{main section}

Introduce a function 
\bean\label{constant}
\lefteqn{
C(z,\la)= \frac{[k]_q\,!\;[k+1]_q\,!}{p^{k+1}\; (k+1)!}\; \Gamma(k/p)\;
\exp(\mu\sum_{a=1}^{k+1}t_a/p)\times} 
\\&&
\times
\prod_{1\leq a< b-1\leq k}\frac{\Gamma((t_a-t_b-1)/p)}{\Gamma((t_a-t_b+1)/p)}
\,\prod_{a=2}^{k+1}
\frac{\Gamma((t_a-z_1+\la_1)/p)}{\Gamma((t_a-z_1-\la_1)/p)}
\,\prod_{a=2}^n\prod_{b=1}^{k+1}
\frac{\Gamma((t_b-z_a+\la_a)/p)}{\Gamma((t_b-z_a-\la_a)/p)},\notag
\eean
where we set $t_a=z_1+\la_1-a+1$, $a=1,\dots, k+1$, cf. formula (16) in \cite{MV}.

\begin{thm}\label{la1=k} Let $z,\la\in\C^n$ and $2\la_1=k\in \Z_{\geq 0}$. Assume $z_a-z_b+\la_a+\la_b\not\in\Z_{\geq 0}\oplus p\Z$ for all $a>b$, $a,b=1,\dots,n$. Let $q=\exp(\pi i/p)$ and $p\not \in \Q$.
Let $\la_1^\prime=-\la_1-1$, and $\tilde{z}=(z_2,\dots,z_n),\;\tilde{\la}=(\la_2,\dots,\la_n)$.
Then there exist a map 
$\alpha(z,\la): L_q(z_1,\la_1)\T V_q(\tilde{z},\tilde{\la})\to L(z_1,\la_1)\T V(\tilde{z},\tilde{\la})$
such that the following diagram is commutative and its columns and rows form exact sequences:
\bean\label{diagram la=k}
\renewcommand{\arraystretch}{1.5}
\begin{array}{ccccccc}
&&{0\qquad\qquad}&&{0\qquad\qquad}&&
\\
&&{\uparrow \qquad\qquad} && {\uparrow \qquad\qquad}&&
\\
0&{\longrightarrow}&{L_q(z_1,\la_1)\T V_q(\tilde{z},\tilde{\la})}& 
{\stackrel{\alpha(z,\la)}{\longrightarrow}}
&{L(z_1,\la_1)\T V(\tilde{z},\tilde{\la})}&{\longrightarrow}&0
\\
&&{\qquad \uparrow Sh_q\T \Id\qquad} && {\qquad \uparrow Sh\T\Id\;\qquad}&&
\\
0&{\longrightarrow}&{V_q(z_1,\la_1)\T V_q(\tilde{z},\tilde{\la})}&
{\stackrel{{\rm qKZ}(z,\la)}{\longrightarrow}}&
{V(z_1,\la_1)\T V(\tilde{z},\tilde{\la})}&{\longrightarrow}&0
\\
&&{\qquad\uparrow\; \iota_q(z,\la)\T \Id \;} && {\qquad\uparrow \; \iota(z,\la)\T \Id\;\; }&&
\\
0&{\longrightarrow}&{V_q(z_1,\la_1^\prime)\T V_q(\tilde{z},\tilde{\la})}&
{\stackrel{C(z,\la)\,{\rm qKZ}(z,\la^\prime)}{\longrightarrow}}&
{V(z_1,\la_1^\prime)\T V(\tilde{z},\tilde{\la})}&{\longrightarrow}&0,
\\
&&{\uparrow \qquad\qquad} && {\uparrow \qquad\qquad}&&
\\
&&{0\qquad\qquad}&&{0\qquad\qquad}&&
\end{array}
\eean
where $C(z,\la)$ is given by \Ref{constant}.
\end{thm}
\begin{proof}
Consider the hypergeometric integral
\be
I(w,W_{\bar{l}})(z,\la)= B^q_{\bar{l}}(\la)\;\int\limits_{{\mathcal K}(z,\la;p)}\;\Phi_l(t,z,\la)w(t,z,\la)W_{\bar{l}}(t,z,\la)\,d^lt, \qquad l_1>k,
\ee
at $\la_1=k/2$. The function $B^q_{\bar{l}}(\la)$ as a function of $\la_1$ has a zero of the first order. The integral $\int\;\Phi_l(t,z,\la)w(t,z,\la)W_{\bar{l}}(t,z,\la)\,d^lt$ has a pole of the first order and the residue is
\bea
\lefteqn{\res_{2\la_1=k}
\int\limits_{{\mathcal K}(z,\la;p)}\;\Phi_l(t,z,\la)w(t,z,\la)W_{\bar{l}}(t,z,\la)\,d^lt={l\choose k+1}\times}
\\&&
\times
\int \limits_{{\mathcal K}(z,\la^\prime;p)}
\res_{t_{k+1}=z_1-\la_1}\ldots
\res_{t_2=z_1+\la_1-1}
\res_{t_1=z_1+\la_1} \Phi_l(t,z,\la)w(t,z,\la)W_{\bar{l}}(t,z,\la) d^{l-k-1}t,
\eea
where $\la^\prime=(\la_1^\prime,\la_2,\dots,\la_n)$, see the proof of Theorem 10 in \cite{MV}. We have
\bea
\lefteqn{
\res_{t_{k+1}=z_1-\la_1}\ldots
\res_{t_2=z_1+\la_1-1}
\res_{t_1=z_1+\la_1}\; \frac{{l\choose k+1}\,\Phi_l(t,z,\la)}{(t-z_1-\la_1)\,\sin(\pi(t-z_1-\l_1)/p)}=}
\\&&
=C(z,\la) \,\Phi_{l-k-1}((t_{k+2},\dots,t_l),z,\la^\prime)\,
\prod_{a=1}^{k+1}\prod_{b=k+2}^l\;\frac{t_a-t_b+1}{t_a-t_b-1}\;\frac{\sin(\pi(t_a-t_b+1)/p)}{\sin(\pi(t_a-t_b-1)/p)}.
\eea

The commutativity of the bottom square of diagram \ref{diagram la=k} now follows from Theorems \ref{iota} and \ref{q iota}. 

There exist the unique map $L_q(z_1,\la_1)\T V_q(\tilde{z},\tilde{\la}) \to L(z_1,\la_1)\T V(\tilde{z},\tilde{\la})$ such that the digram \ref{diagram la=k} is commutative. 
\end{proof}


\begin{remark}
A version of Theorem \ref{la1=k} was implicitely used in \cite{MV} to construct solutions of the qKZ equation with values in tensor products of finite dimensional representations $L(\la_1)\T\dots\T L(\la_n)$, $2\la\in\Z^n_{\geq 0}$.
\end{remark}

\begin{corollary}\label{qKZ same kind}
Let $z,\la\in\C^n$ be in the first trigonometric resonance and in the first rational resonance. Then the map 
${\rm qKZ}(z,\la): \; V_q(z,\la)\to V(z,\la)$ is a well defined isomorphism of linear spaces. Moreover, it maps the submodule $V_q(z^\prime,\la^\prime)\subset V_q(z,\la)$ onto the submodule $V(z^\prime,\la^\prime)\subset V(z,\la)$. The map ${\rm qKZ}(z,\la)$ restricted to the submodule $V_q(z^\prime,\la^\prime)$ coincides with the map 
${\rm qKZ}(z^\prime,\la^\prime)$ up to a non-zero scalar multiplier depending on the choice of the inclusions $V_q(z^\prime,\la^\prime)\hookrightarrow V_q(z,\la)$ and
$V(z^\prime,\la^\prime)\hookrightarrow V(z,\la)$.
\end{corollary}
\begin{proof}
Corollary \ref{qKZ same kind} follows from Theorem \ref{la1=k} and Lemmas \ref{same kind}, \ref{q same kind}.
\end{proof}

\begin{corollary}\label{qKZ 1 to 2}
Let $z,\la\in\C^n$ be in the first trigonometric resonance and in the second rational resonance. Then the map 
${\rm qKZ}(z,\la): \; V_q(z,\la)\to V(z,\la)$ is a well defined linear map. 
The kernel of this map is 
the submodule $V_q(z^\prime,\la^\prime)\subset V_q(z,\la)$ and the image is the proper submodule in $V(z,\la)$.
\end{corollary}
\begin{proof}
Corollary \ref{qKZ 1 to 2} follows from  Theorem \ref{la1=k} and Lemmas \ref{1 to 2}, \ref{q 1 to 2}.
\end{proof}

\begin{corollary}\label{qKZ 2 to 1}
Let $z,\la\in\C^n$ be in the second trigonometric resonance and in the first rational resonance. Then the map 
${\rm qKZ}(\tilde{z},\tilde{\la}): \; V_q(\tilde{z},\tilde{\la})\to V(\tilde{z},\tilde{\la})$ has a simple pole at $\tilde{z}=z,\tilde{\la}=\la$. Let $Res$ be the residue  of the map ${\rm qKZ}(\tilde{z},\tilde{\la})$ at $\tilde{z}=z,\tilde{\la}=\la$. The kernel of the map $Res$ is 
the nontrivial submodule $U_q$ of $V_q(z,\la)$ and the image is the submodule 
$V(z^\prime,\la^\prime)\subset V(z,\la)$. Thus, up to a scalar multiplier, 
depending on the choice of the factor map $V_q(z,\la)\to V_q(z^\prime,\la^\prime)\simeq V_q(z,\la)/U_q$, and the inclusion $V(z^\prime,\la^\prime)\hookrightarrow V(z,\la)$,
 the map $Res$ defines a homomorphism $V_q(z^\prime,\la^\prime)\to V(z^\prime,\la^\prime)$. The scalar multiplier can be chosen so that the map $Res$ coincides with the isomorphism
${\rm qKZ}(z^\prime,\la^\prime)$.
\end{corollary}
\begin{proof}
Corollary \ref{qKZ 2 to 1} follows from Theorem \ref{la1=k} and Lemmas \ref{2 to 1}, \ref{q 2 to 1}. 
\end{proof}

\begin{example}\label{z1-z2=k}
\rm
Let $z,\la\in\C^n$ be in the second trigonometric resonance and in the first rational resonance. Let $(\la^{(1,2)}_{\;id})_1=z_1-z_2+\la_1+\la_2=k\in \Z_{\geq 0}$. Denote $u=z^{(1,2)}_{\;id}, \om=\la^{(1,2)}_{\;id}$, $\tilde{z}=(z_3,\dots,z_n)$, $\tilde{\la}=(\la_3,\dots,\la_n)$. Then there exists the commutative diagram
\bea
\renewcommand{\arraystretch}{1.5}
\begin{array}{ccccccc}
&&{0\qquad}&&{0\qquad}&&
\\
&&{\downarrow \qquad} && {\uparrow \qquad}&&
\\
&&{L_q(u_1,\om_1)\T V_q(u_2,\om_2)\T V_q(\tilde{z},\tilde{\la})}& 
{\stackrel{0}{\longrightarrow}}
&{L(u_1,\om_1)\T V(u_2,\om_2)\T V(\tilde{z},\tilde{\la})}&&
\\
&&{\downarrow\qquad } && {\uparrow \qquad}&&
\\
&&{V_q(z_1,z_2,\la_1,\la_2)\T V_q(\tilde{z},\tilde{\la})}&
{\stackrel{Res}{\longrightarrow}}&
{V(z_1,z_2,\la_1,\la_2)\T V(\tilde{z},\tilde{\la})}&&
\\
&&{\downarrow\; \qquad } && {\uparrow \qquad }&&
\\
0&{\to}&{V_q(u_1,\om_1^\prime)\T V_q(u_2,\om_2)\T V_q(\tilde{z},\tilde{\la})}&
{\stackrel{{\rm qKZ}}{\longrightarrow}}&
{V(u_1,\om_1^\prime)\T V(u_2,\om_2)\T V(\tilde{z},\tilde{\la})}&{\to}&0.
\\
&&{\downarrow \qquad} && {\uparrow \qquad}&&
\\
&&{0\qquad}&&{0\qquad}&&
\end{array}
\eea
Here the left column is an exact short sequence of $\Ugh$ modules and the right column is an exact short sequence of Yangian modules, 
$Res=\res_{z_1=z_2-\la_1-\la_2+k} \;{\rm qKZ}(\tilde{z},\tilde{\la})$.
\end{example}

\subsection{Singularities of hypergeometric solutions}\label{singularities section}

\begin{thm}\label{singularities}
Let $p$ be a complex number such that $\Real p <0$ and $p\not\in\Q$. 
The map qKZ$(z,\la): V_q(z,\la)\to V(z,\la)$ is a well defined isomorhism of vector spaces for all $z,\la$ except for the following two cases.

i) The map qKZ$(z,\la)$ has a nontrivial simple pole at the hyperplanes $z_i-z_j+\la_i+\la_j=m-ps$, $m,s\in\Z_{\geq 0}$, $i<j$.

ii) The map qKZ$(z,\la)$ has a nontrivial kernel at the hyperplanes $z_i-z_j+\la_i+\la_j=m+ps$, $m,s\in\Z_{\geq 0}$, $j<i$.
\end{thm}
\begin{proof}
By Corollary 5.25 in \cite{TV} the map qKZ$(z,\la)$ is a linear isomorphism if $z_i-z_j+\la_i+\la_j\not\in\Z_{\geq 0}+p\Z$. 

Consider the case $i=1$, $j=2$.
By Corollary \ref{qKZ 2 to 1}, the map qKZ$(z,\la)$ has a pole at the hyperplane $z_1-z_2+\la_1+\la_2=k\in\Z_{\geq 0}$. We have
\be
(H_1(z,\la))^s{\rm qKZ}(z,\la)={\rm qKZ}(z_1+sp,z_2,\dots,z_n,\la)), \qquad s\in\Z,
\ee
where $H_1$ is the first qKZ operator given by \Ref{qKZ operators}. The operator $H_1(z,\la)$ has a non-trivial kernel at the hyperplane $z_1-z_2+\la_1+\la_2=k\in\Z_{\geq 0}$. It is easy to see that at this hyperplane for generic $z_3,\dots, z_n$, the product $H_1(z,\la)qKZ(z,\la)$ is a well defined nondegenerate operator. 

Note that for generic $z_3,\dots ,z_n$, the operator $H_1(z,\la)$ is an isomorphism if $z_1-z_2+\la_1+\la_2=k+sp$, $s\neq 0, z\in\Z$. Hence, the map qKZ$(z,\la)$ is an isomorpism at $z_1-z_2+\la_1+\la_2\in\Z+p\Z_{>0}$ and has a nontrivial pole at $z_1-z_2+\la_1+\la_2\in\Z+p\Z_{\leq 0}$.

The case of generic $i,j\in \{1,\dots, n\}$ is done similarly.
\end{proof}

\begin{remark}
In fact Theorem \ref{singularities} combined with the results of Section \ref{main section} allows to describe all the singularities of the hypergeometric solutions. Indeed, consider a hypergeometric solution $\Psi_{\bar{l}}(z,\la)$ given by \Ref{Psi}. It has poles of the first order at the hyperplanes $z_i-z_j+\la_i+\la_j=m-ps$, $s\in\Z_{\geq 0}$, $m=0,1\dots, l-1$, $i<j$. 

Consider the hyperplane $z_1-z_2+\la_1+\la_2=m$. The residue $\res_{z_1=z_2-\la_1-\la_2+k} \Psi_{\bar{l}}(z,\la)$ is a function with values in the proper Yangian submodule $V(u,\om)\subset V(z,\la)$, cf. Example \ref{z1-z2=k}. In fact this function is a linear combination of hypergeometric solutions $\Psi_{\bar{m}}(u,\om)$ with $m=l-k-1$ integrations. This linear combination is determined by the image of the vector 
$f^{l_1}v_1\T \dots\T f^{l_n}v_n\in V_q(z,\la)$ under the factorization map of $\Ugh$ modules $V_q(z,\la)\to V_q(u,\om)$.

The residue of the hypergeometric function $\Psi_{\bar{l}}(z,\la)$ at the hyperplanes  $z_1-z_2+\la_1+\la_2=m-ps$, $s=1,2,\dots$, is computed by applying $s$ times the qKZ operator $H_2$ to the residue of $\Psi_{\bar{l}}(z,\la)$ at 
the hyperplane $z_1-z_2+\la_1+\la_2=m$.
\end{remark}

\begin{remark}\label{kappa=1}
In this paper we treated the case $|\kappa|\neq 1$. The case $\kappa=1$ is important for applications and is done in a similar way. In this case the qKZ map is defined on the subspaces of singular vectors,
\be
{\rm qKZ}(z,\la;\,p,\kappa=1): (V^q(z,\la))^{sing}\to (V(z,\la))^{sing},
\ee
where $(V^q(z,\la))^{sing}=\Ker\; e_q\subset V_q(z,\la)$, $(V(z,\la))^{sing}=\Ker\;e\subset V(z,\la)$. We get the same statements as in Corollaries \ref{qKZ same kind}-\ref{qKZ 2 to 1}, where all the tensor products of modules are replaced with the corresponding subspaces of singular vectors. 

However, in the case $\kappa=1$, the map qKZ$(z,\la)$ has nontrivial degeneracies which do not come from singularities of the qKZ equation, see \cite{MV2},\cite{MV3}. 
\end{remark}

\begin{remark}\label{finite case remark}
One can prove a statement similar to Corollary \ref{qKZ 2 to 1} for the case $2\la\in\Z^n_{\geq 0}$ with Verma modules $V(\la_i), V_q(\la_i)$ replaced by finite dimensional modules $L(\la_i), L_q(\la_i)$. In this case, it is sufficient to assume that $q$ is not a root of unity of a small order. Namely, it is sufficient to assume $q^a\neq 1$ for $a=1,\dots,\max(2\la_1,\dots,2\la_n, k, 2\la_1+2\la_2-k)$. 
\end{remark}

\bigskip

{\it Mathematical Sciences Research Institute, 1000 Centennial Drive, Berkeley, CA 94720-5070}

{\it Department of Mathematics, University of North Carolina at Chapel Hill, Chapel Hill, NC 27599-3250, USA.}

{\it E-mail addresses:} {\rm mukhin@@msri.org,
av@@math.unc.edu}

\end {document}